\documentclass[12pt,a4paper]{amsart}
\def\amsArticle

\usepackage{t1enc}
\usepackage[latin1]{inputenc}
\usepackage[english]{babel}
\usepackage[linktocpage=true]{hyperref}
\usepackage{graphicx}
\usepackage{latexsym}
\usepackage{mathtools}
\usepackage{amsmath,amsfonts,amssymb,amsthm}
\usepackage{enumitem}
\ifdef{\amsArticle}{
}

\everymath{\displaystyle}
\begin{document}

\title{Lipshitzian Vector Fields, Upper Gradients And Distributional Derivatives}
\author{Sergio Venturini}

\ifdef{\amsArticle} {
\address{
	Sergio Venturini,
	Dipartimento di Matematica,
	Universit\`{a} di Bologna,
	Piazza di Porta S. Donato 5 ---I-40127 Bologna,
	Italy
}
\email{sergio.venturini@unibo.it}
}

\ifdef{\elseArticle}{
\address{
	Dipartimento di Matematica,
	Universit\`{a} di Bologna,
	
	Piazza di Porta S. Donato 5 ---I-40127 Bologna,
	Italy
	
	sergio.venturini@unibo.it}
}

\date{\today}
%\maketitle

%/intest
%\magnification=\magstep1
%\font\chapt=cmbx10 scaled 1700
%\font\sc=cmcsc10
%\font\ninerm=cmr8
%\font\nineit=cmmi8
%\font\ninesl=cmsl8
%\font\titolo=cmbx9 scaled \magstep 3
%
%
%numeri....
\def\R{{\rm I\kern-.185em R}}
\def\RR{\mathbb{R}}
\def\C{{\rm\kern.37em\vrule height1.4ex width.05em depth-.011em\kern-.37em C}}
\def\CC{\mathbb{C}}
\def\N{{\rm I\kern-.185em N}}
\def\NN{\mathbb{N}}
\def\Z{{\bf Z}}
\def\ZZ{\mathbb{Z}}
\def\Q{\mathbb{Q}}
\def\P{{\rm I\kern-.185em P}}
\def\H{{\rm I\kern-.185em H}}
%
%insiemi
\def\Aleph{\aleph_0}
\def\ALEPH#1{\aleph_{#1}}
\def\sset{\subset}\def\ssset{\sset\sset}
%
%funzioni
\def\bar#1{\overline{#1}}
\def\dim{\mathop{\rm dim}\nolimits}
\def\half{\textstyle{1\over2}}
\def\Half{\displaystyle{1\over2}}
\def\mlog{\mathop{\half\log}\nolimits}
\def\Mlog{\mathop{\Half\log}\nolimits}
\def\Det{\mathop{\rm Det}\nolimits}
\def\Hol{\mathop{\rm Hol}\nolimits}
\def\Aut{\mathop{\rm Aut}\nolimits}
\def	\Re{\mathop{\rm Re}\nolimits}
\def\Im{\mathop{\rm Im}\nolimits}
\def\Ker{\mathop{\rm Ker}\nolimits}
\def\Fix{\mathop{\rm Fix}\nolimits}
\def\Exp{\mathop{\rm Exp}\nolimits}
\def\sp{\mathop{\rm sp}\nolimits}
\def\id{\mathop{\rm id}\nolimits}
\def\Rank{\mathop{\rm rk}\nolimits}
\def\Trace{\mathop{\rm Tr}\nolimits}
\def\Res{\mathop{\rm Res}\limits}
\def\Divergence{\mathop{\rm div}\nolimits}
\def\cancel#1#2{\ooalign{$\hfil#1/\hfil$\crcr$#1#2$}}
\def\prevoid{\mathrel{\scriptstyle\bigcirc}}
\def\void{\mathord{\mathpalette\cancel{\mathrel{\scriptstyle\bigcirc}}}}
\def\n{{}|{}\!{}|{}\!{}|{}}
\def\abs#1{\left|#1\right|}
\def\norm#1{\left|\!\left|#1\right|\!\right|}
\def\nnorm#1{\left|\!\left|\!\left|#1\right|\!\right|\!\right|}
\def\Norm#1{\Bigl|\!\Bigl|#1\Bigr|\!\Bigr|}
%
%integrali superiore ed inferiore
\def\upperint{\int^{{\displaystyle{}^*}}}
\def\lowerint{\int_{{\displaystyle{}_*}}}
\def\Upperint#1#2{\int_{#1}^{{\displaystyle{}^*}#2}}
\def\Lowerint#1#2{\int_{{\displaystyle{}_*}#1}^{#2}}
%
%altro
\def\rem #1::#2\par{\medbreak\noindent{\bf #1}\ #2\medbreak}
\def\proclaim #1::#2\par{\removelastskip\medskip\goodbreak{\bf#1:}
\ {\sl#2}\medskip\goodbreak}
\def\ass#1{{\rm(\rmnum#1)}}
\def\assertion #1:{\Acapo\llap{$(\rmnum#1)$}$\,$}
\def\Assertion #1:{\Acapo\llap{(#1)$\,$}}
\def\acapo{\hfill\break\noindent}
\def\Acapo{\hfill\break\indent}
\def\prova{\removelastskip\par\medskip\goodbreak\noindent{\it Dimostrazione.\/\ }}
\def\qed{{$\Box$}\par\smallskip}
\def\BeginItalic#1{\removelastskip\par\medskip\goodbreak
\noindent{\it #1.\/\ }}
\def\iff{if, and only if,\ }
\def\sse{se, e solo se,\ }
\def\rmnum#1{\romannumeral#1{}}
\def\Rmnum#1{\uppercase\expandafter{\romannumeral#1}{}}
\def\smallfrac#1/#2{\leavevmode\kern.1em
\raise.5ex\hbox{\the\scriptfont0 #1}\kern-.1em
/\kern-.15em\lower.25ex\hbox{\the\scriptfont0 #2}}
%
%delimitatori
\def\Left#1{\left#1\left.}
\def\Right#1{\right.^{\llap{\sevenrm
\phantom{*}}}_{\llap{\sevenrm\phantom{*}}}\right#1}
\def\definedby{\mathrel{\mathop:}=}
\def\newpi{{\pi\mskip -7.8 mu \pi}} %Richard Palais newPi
%bibliografia
%\def\dimens{3em}
%\newcount\qqrefno
%\qqrefno=0
%\def\qqrefnum{\global\advance\qqrefno by 1
%\noindent\rlap{[\number\qqrefno]}\hbox to \dimens{}\hangindent=\dimens}
%\def\references{\bigskip\noindent{\bf References.}\bigskip}
%\def\art #1 : #2 ; #3 ; #4 ; #5 ; #6. \par{\qqrefnum
%#1, {\sl#2}, #3, {\bf#4}, (#5), #6.\par\smallskip}
%\def\book #1 : #2 ; #3 ; #4. \par{\qqrefnum#1, {\bf#2}, #3, #4.\par\smallskip}
%\def\freeart #1 : #2 ; #3. \par{\qqrefnum#1, {\sl#2}, #3.\par\smallskip}
%
%bibliografiabis
\def\dimens{3em}
\def\symb[#1]{\noindent\rlap{[#1]}\hbox to \dimens{}\hangindent=\dimens}
\def\references{\bigskip\noindent{\bf References.}\bigskip}
\def\art #1 : #2 ; #3 ; #4 ; #5 ; #6. \par{#1, 
{\sl#2}, #3, {\bf#4}, (#5), #6.\par\smallskip}
\def\book #1 : #2 ; #3 ; #4. \par{#1, {\bf#2}, #3, #4.\par\smallskip}
\def\freeart #1 : #2 ; #3. \par{#1, {\sl#2}, #3.\par\smallskip}
%
%
%bibliografia per riviste tedesche
%\def\dimens{3em}
%\newcount\qqrefno
%\qqrefno=0
%\def\qqrefnum{\global\advance\qqrefno by 1
%\noindent\hbox to \dimens{}\hangindent=\dimens
%\llap{\number\qqrefno.\ \ \ }}
%\def\references{\bigskip\noindent{\bf References.}\bigskip}
%\def\art #1 : #2 ; #3 ; #4 ; #5 ; #6. \par{\qqrefnum
%#1: #2. #3 {\bf#4}, #6, (#5)\par\smallskip}
%\def\book #1 : #2 ; #3. \par{\qqrefnum#1: #2. #3.\par\smallskip}
%\def\freeart #1 : #2 ; #3. \par{\qqrefnum#1: #2. #3.\par\smallskip}
%
%
%nome e indirizzi.
\def\name{\hbox{Sergio Venturini}}
\def\snsaddress{\indent
\vbox{\bigskip\bigskip\bigskip
\name
\hbox{Scuola Normale Superiore}
\hbox{Piazza dei Cavalieri, 7}
\hbox{56126 Pisa (ITALY)}
\hbox{FAX 050/563513}}}
\def\cassinoaddress{\indent
\vbox{\bigskip\bigskip\bigskip
\name
\hbox{Universit\`a di Cassino}
\hbox{via Zamosch 43}
\hbox{03043 Cassino (FR)}
\hbox{ITALY}}}
\def\bolognaaddress{\indent
\vbox{\bigskip\bigskip\bigskip
\name
\hbox{Dipartimento di Matematica}
\hbox{Universit\`a di Bologna}
\hbox{Piazza di Porta S. Donato 5}
\hbox{40127 Bologna (BO)}
\hbox{ITALY}
\hbox{sergio.venturini@unibo.it}
}}
\def\homeaddress{\indent
\vbox{\bigskip\bigskip\bigskip
\name
\hbox{via Garibaldi, 7}
\hbox{56124 Pisa (ITALY)}}}
\def\doubleaddress{
\vbox{
\hbox{\name}
\hbox{Universit\`a di Cassino}
\hbox{via Zamosch 43}
\hbox{03043 Cassino (FR)}
\hbox{ITALY}
\smallskip
\hbox{and}
\smallskip
\hbox{Scuola Normale Superiore}
\hbox{Piazza dei Cavalieri, 7}
\hbox{56126 Pisa (ITALY)}
\hbox{FAX 050/563513}}}
\def\sergio{{\rm\bigskip
\centerline{Sergio Venturini}
\leftline{\bolognaaddress}
\bigskip}}
%
%
%allineamenti formule matematiche
%
%\centeredeq{...}
%                      |*  [ ]&*   [ ]&....|
%                      ..................
%\leftcenteredeq
%|*  [ ]&*   [ ]&....|
%..................
%
%\def\centeredeq#1{\vcenter{\halign
%{&$\strut\displaystyle{##}\hfil\strut\qquad$\cr#1}}
%\def\leftcenteredeq#1{\halign
%{&$\strut\displaystyle{##}\hfil\strut\qquad$\cr#1}}
%
%@@@@@@@
%
%alfabeto greco.
\def\a{\alpha}
\def\bg{\beta}
\def\g{\gamma}
\def\G{\Gamma}
\def\dg{\delta}
\def\D{\Delta}
\def\e{\varepsilon}
\def\eps{\epsilon}
\def\z{\zeta}
\def\th{\theta}
\def\T{\Theta}
\def\k{\kappa}
\def\lg{\lambda}
\def\Lg{\Lambda}
\def\m{\mu}
\def\n{\nu}
\def\r{\rho}
\def\s{\sigma}
\def\Sg{\Sigma}
\def\ph{\varphi}
\def\Ph{\Phi}
\def\x{\xi}
\def\om{\omega}
\def\Om{\Omega}

%macrosl

%\theoremstyle{plain}% default
%\newtheorem{thm}{Theorem}[section]
%\newtheorem{lem}[thm]{Lemma}
%\newtheorem{prop}[thm]{Proposition}
%\newtheorem*{cor}{Corollary}
%\newtheorem*{KL}{Klein's Lemma}

%\theoremstyle{definition}
%\newtheorem{defn}{Definition}[section]
%\newtheorem{conj}{Conjecture}[section]
%\newtheorem{exmp}{Example}[section]

%\theoremstyle{remark}
%\newtheorem*{rem}{Remark}
%\newtheorem*{note}{Note}
%\newtheorem{case}{Case}

%\theoremstyle{plain} %default

%english
\newtheorem{theorem}{Theorem}[section]
\newtheorem{proposition}[theorem]{Proposition}
\newtheorem{lemma}[theorem]{Lemma}
\newtheorem{corollary}[theorem]{Corollary}

%italiano
\newtheorem{teorema}{Teorema}[section]
\newtheorem{proposizione}[teorema]{Proposizione}
\newtheorem{corollario}[teorema]{Corollario}

\newtheorem{definition}[theorem]{Definition}

\newtheorem{definizione}[teorema]{Definizione}

\newtheorem{remark}[theorem]{Remark}

\newtheorem{osservazione}[teorema]{Osservazione}
\newtheorem{esempio}[teorema]{Esempio}
\newtheorem{esercizio}[teorema]{Esercizio}
\newtheorem{congettura}[teorema]{Congettura}

% !TeX encoding = UTF-8

\def\Dim{n}
\def\ManBase{M}
\def\VField{X}
\def\VFieldA{X}
\def\VFieldB{Y}
\def\VFlowBase{\gamma}
\def\VFlowMap#1#2{\VFlowBase_{#2}^{#1}}
\def\VFlow#1#2#3{\VFlowMap{#1}{#2}(#3)}
\def\DFlow{\Gamma}
\def\DFlowX#1{\DFlow^{#1}}
\def\DFlowXPlus#1{\DFlow_{+}^{#1}}
\def\DFlowP#1#2{{\DFlow}_{#1}^{#2}}
\def\pMan{p}
\def\pManA{p}
\def\pManB{q}
\def\FlowVar{t}
\def\FlowVarA{t}
\def\FlowVarB{s}
\def\FlowVarI{u}
\def\FlowVarM{h}

\def\funcTest{u}
\def\funcA{f}
\def\funcB{g}
\def\UGradA{h}
\def\LipC#1{{\rm Lip}_{0}(#1)}
\def\LOneLoc#1{L^1_{\rm loc}({#1})}
\def\VFSubSpace{S}

\def\MDer#1#2#3{\varphi_{#1}\left({#2},{#3}\right)}
\def\DFGenA{\beta}

\def\CutOffA{\rho}
\def\Cutoff{\rho}
\def\VFieldLoc{\tilde{\VFieldE}}
\def\OSubsetRK{U}
\def\VFieldCS{\tilde\VFieldE}

\def\VFieldE{\VField}
\def\EucNorm#1{\abs{{#1}}}
\def\VFComp{a}
\def\LipVF{L}
\def\OpenE{\Omega}
\def\SubOpenE{\OpenE'}
\def\pEuc{x}
\def\pEucA{x}
\def\pEucB{y}
\def\VFAdvance{\lambda}
\def\EucBall#1#2{B({#1},{#2})}
\def\EucBallB#1#2{B\bigl({#1},{#2}\bigr)}

\def\MeanOp#1#2#3{\textrm{M}_{#2}^{#1}#3}
\def\DeltaR#1#2#3{\Delta_{#2}^{#1}#3}
\def\VFlowW#1#2{\gamma_{#2}^{#1}}
\def\VFlow#1#2#3{\gamma_{#2}^{#1}(#3)}
\def\VFlowB#1#2#3{\gamma_{#2}^{#1}\bigl(#3\bigr)}
\def\VFJac#1#2#3{J_{#2}^{#1}(#3)}
\def\VFJacB#1#2#3{J_{#2}^{#1}\bigl(#3\bigr)}
\def\TestFEuc{u}
\def\COneFTest{u}
\def\GPFlow#1#2{T^{#1}_{#2}}
\def\GFZ{Z}
\def\Dom{{\rm Dom}}

\def\VFIMin{1}
\def\VFIMax{k}
\def\VFIndex{j}
\def\VConstX{c}

\def\MatrSpace{M_\Dim(\RR)}
\def\MatrNorm#1{\norm{#1}}
\def\MatrZero{0_\Dim}
\def\MatrId{I_\Dim}

\def\LimVFZ{\lim_{\FlowVar\to0}}

\def\cbLipz{Lipshitz }
\def\cbLipn{Lipshitzian }

\def\BanachS{B}
\def\BNorm#1{\norm{#1}}
\def\BElem{f}
\def\BElemZ{g}
\def\BStarElem{\varphi}
\def\BSMap{T}
\def\BGen{Z}
\def\BDomGen#1{\textrm{\rm Dom({#1})}}
\def\SCore{D}
\def\wstar{\({}^*\ \)}
\def\weakstar{weak\wstar}
\def\BDuality#1#2{\langle{#1},{#2}\rangle}

\def\ChParam{\alpha}
\def\cbDot#1{{#1}'}
\def\xVar{\FlowVar}

\def\Int{\int_\OpenE}
\def\NullSet{N}

%%%%%%%%%%%%%%

\def\Case#1{\((\rmnum{#1})\)}

\def\stmJBouned{1}
\def\stmJwStarToDiv{2}

\def\stmGFXLoc{1}
\def\stmGFLLoc{2}
\def\stmGFZLoc{3}

\def\stmGFX{1}
\def\stmGFL{2}
\def\stmGFZ{3}

\def\cbRomanLabel{\((\)\roman*\()\)}
{%begin unit

}%end unit
\def\cbFullBase{true}

\def\cbVThree{defined}
%
%\def\cbAppendixA{defined}
%
%
% !TeX encoding = UTF-8
% !TeX spellcheck = en_GB

%\section{\label{section:Intro}Introduction}
% !TeX encoding = UTF-8
% !TeX spellcheck = en_GB

{%begin unit
\begin{abstract}
%In this paper we study the relationships
%between
%the distributional derivatives and the upper gradients
%witrh respect to a locally Lipshitzian vector fields on open set of the Euclidean space.
We prove that given a locally integrable function \(\funcA\)
on an open set of an Euclidean space
%and a locally Lipshitzian vector fields \(\VFieldE\) 
the distributional derivative \(\VFieldE\funcA\) %of \(\funcA\)
with respect to a locally Lipshitzian vector field \(\VFieldE\)
is locally integrable if, and only if,
the function \(\funcA\) admits a locally integrable upper gradient
along the vector field \(\VFieldE\);
in this case
\(\VFieldE\funcA\) coincides with the Lie derivative \(L_\VFieldE\funcA\)
and \(\abs{\VFieldE\funcA}\)
%the modulus of the distributional derivative
is the least upper gradient of the function \(\funcA\).
%No kind of continuity is required on the function \(\funcA\). 
Applications to systems of locally Lipshitzian vector fields are given.
\end{abstract}
}%end unit

\maketitle

\tableofcontents

\section{\label{section:Main}Introduction}
%\subsection{Motivations}
%\input{src/LipVFieldsMotivations}
% !TeX encoding = UTF-8
% !TeX spellcheck = en_GB
{%begin unit
%I risultati principali di questo lavoro si inseriscono nella cerchia di idee riguardanti le funzioni localmente integrabili che apparterngono agli spazi di Sobolev associati a campi di vettori localmente Lipshitziani su spazi Euclidei.
%I risultati principali di questo lavoro si inseriscono nella cerchia di idee riguardanti le funzioni localmente integrabili che apparterngono agli spazi di Sobolev definiti mediante campi di vettori localmente Lipshitziani su spazi Euclidei.
%Si caratterizzano tali funzioni mediante "uper gradients" lungo campi vettoriali definiti in modo simile a quelli 
The main results of this work fit into the circle of ideas regarding
the analysis on Carnot-Carath\'eodory spaces from the point of view of upper gradients
and the
(locally) integrable functions that belong to
Sobolev spaces associated to locally Lipshitzian vector fields on Euclidean spaces.
	}%end unit
\subsection{Statements of main results}
%\ifdef{\cbVOne}{
%	\input{src/LipVFieldsIMain}
%}{}
%\ifdef{\cbVTwo}{
%	\input{src/LipVFieldsIMainVTwo}
%}{}
%\ifdef{\cbVThree}{
% !TeX encoding = UTF-8
% !TeX spellcheck = en_GB
{%begin unit

A vector field on an open set \(\OpenE\sset\RR^\Dim\)
may be identified with a map
\(\VFieldE:\OpenE\to\RR^\Dim\).

If \(\VFComp_1,\ldots,\VFComp_\Dim\) are the components of \(\VFieldE\)
we also write
\begin{equation*}
	\VFieldE=\VFComp_1\dfrac{\partial}{\partial\pEuc_1}
	+\cdots+\VFComp_\Dim\dfrac{\partial}{\partial\pEuc_\Dim}
	.
\end{equation*}
The action of the vector field \(\VFieldE\) on a function \(\funcA\in C^1(\OpenE)\)
is given by
\begin{equation*}
	\VFieldE\funcA(\pEuc)=\VFComp_1(\pEuc)\dfrac{\partial\funcA(\pEuc)}{\partial\pEuc_1}
	+\cdots+\VFComp_\Dim(\pEuc)\dfrac{\partial\funcA(\pEuc)}{\partial\pEuc_\Dim},
\end{equation*}
%By the classical Rademacher theorem the same formula defines
%a (locally) bounded Borel function when \(\funcA\) is
%a (locally) Lipshitzian function.
and a \(C^1\) curve
\(%\begin{equation*}
	\RR\supset]a,b[\ni\FlowVar\mapsto\gamma(t)\in\OpenE\
\)%\end{equation*}
is an \emph{integral curve} of the vector field \(\VField\)
if satisfies
\(%\begin{equation*}
	\dot\gamma(\FlowVar)=\VField\bigl(\gamma(\FlowVar))\bigr)
\) %\end{equation*}
for each \(\FlowVar\in]a,b[\).

%We note that if \(\funcA\in C^1(\OpenE)\) then
%\begin{equation*}
%	\VFieldE\funcA(\pEuc)=\VFComp_1(\pEuc)\dfrac{\partial\funcA(\pEuc)}{\partial\pEuc_1}
%	+\cdots+\VFComp_\Dim(\pEuc)\dfrac{\partial\funcA(\pEuc)}{\partial\pEuc_\Dim}.
%\end{equation*}
%
If the vector field \(\VFieldE\) is locally
Lipshitzian on the open subset \(\OpenE\subset\RR^\Dim\)
then for each \(\pEuc\in\OpenE\) there exists a unique (maximal) integral curve
%denote by
\begin{equation*}
	\RR\supset\DFlowP\pEuc\VField\ni\FlowVar\mapsto\VFlow{\VField}{\FlowVar}{\pEuc}\in\OpenE
\end{equation*}
%the (unique) maximal integral curve
of the vector field \(\VField\)
that passes through the point \(\pEuc\in\OpenE\) when \(\FlowVar=0\)
%is the unique curve
%which satisfies
%\begin{eqnarray*}
%	&&\dotVFlow{\VField}{\FlowVar}{\pEuc}=
%	\VField\left(\VFlow{\VField}{\FlowVar}{\pEuc}\right),\\
%	&&\VFlow{\VField}{0}{\pEuc}=\pEuc.
%\end{eqnarray*}
and the \emph{flow} associated to the vector field \(\VField\)
is the map
\begin{equation*}
	\RR\times\OpenE\supset\DFlowX\VField\ni(\FlowVar,\pEuc)\mapsto\VFlow{\VField}{\FlowVar}{\pEuc}\in\OpenE;
\end{equation*}
where
\begin{equation*}
	\DFlowX\VField\definedby\bigcup_{\pEuc\in\OpenE}\DFlowP\pEuc\VField\times\{\pEuc\}.
\end{equation*}
It turns out that \(\DFlowX\VField\) is an open subset of \(\RR\times\OpenE\)
containing \(\{0\}\times\OpenE\)
and hence, for each relatively compact open set \(\OpenE'\sset\OpenE\)
there exists \(\FlowVar_0>0\) such that
\(%\)\begin{equation*}
	[-\FlowVar_0,\FlowVar_0]\times\OpenE'\subset\DFlowX\VField.
\)%\end{equation*}

Throughout this paper we will use the notation
\begin{equation*}
	\DFlowXPlus\VField\definedby\bigl\{(\FlowVar,\pEuc)\in\DFlowX\VField\mid\FlowVar>0\bigr\}
\end{equation*}
and ``almost all \((\FlowVar,\pEuc)\in\DFlowXPlus\VField\)''
stands for almost all \((\FlowVar,\pEuc)\in\DFlowXPlus\VField\)
with respect to the Lebesgue meaasure induced on \(\DFlowXPlus\VField\)
as open subset of \(]0,+\infty[\times\RR^\Dim\subset\RR^{\Dim+1}\).

The \emph{difference quotient} associated to the function \(\funcA\in\LOneLoc{\OpenE}\)
with respect to the vector field \(\VField\) is defined
for each \((\FlowVar,\pEuc)\in\DFlowX\VField\), \(\FlowVar\neq0\) as
\begin{equation*}
	\DeltaR\VField\FlowVar\funcA(\pEuc)
	\definedby
	\dfrac{\funcA\bigl(\VFlow{\VField}{\FlowVar}{\pEuc}\bigr)-\funcA(\pEuc)}{\FlowVar}
\end{equation*}
and we also define the \emph{mean operator} by the formula
\begin{equation*}%\label{stm::MeanOpDef}
	\MeanOp{\VFieldE}{\FlowVar}{\funcA}(\pEuc)=\dfrac{1}{\FlowVar}
	\int_0^\FlowVar\funcA\bigl(\VFlow{\VFieldE}{\FlowVarB}{\pEuc}\bigr)\,d\FlowVarB.
\end{equation*}

For convenience we set \(\DeltaR\VField\FlowVar\funcA(\pEuc)=\MeanOp{\VFieldE}{\FlowVar}{\funcA}(\pEuc)=0\)
if \((\FlowVar,\pEuc)\in\RR^{\Dim+1}\setminus\OpenE\).
We also set \(\MeanOp{\VFieldE}{0}{\funcA}=\funcA\).

Inspired by \cite{article:HajlaszKoskela:SobolevMetPoincare},
\cite{article:HeinonenKoskela98},\cite{article:Cheeger99},
we say that a nonnegative function \(\UGradA\in\LOneLoc{\OpenE}\)
is an \emph{upper gradient} for the function \(\funcA\in\LOneLoc{\OpenE}\)
along the vector field \(\VFieldE\)
if
%whenever \([0,\FlowVar_0]\times\OpenE'\subset\DFlowX\VField\)
%and \(0<\FlowVar<\FlowVar_0\)
%then %for almost all \(\pEuc\in\OpenE'\)
\begin{equation*}
	\abs{\funcA\bigl(\VFlow{\VField}{\FlowVar}{\pEuc}\bigr)-\funcA(\pEuc)}
	\leq
	\int_0^\FlowVar\UGradA\bigl(\VFlow{\VFieldE}{\FlowVarB}{\pEuc}\bigr)\,d\FlowVarB,
\end{equation*}
or equivalently
\begin{equation*}
	\abs{\DeltaR{\VFieldE}{\FlowVar}{\funcA}(\pEuc)}
	\leq
	\MeanOp{\VFieldE}{\FlowVar}{\UGradA}(\pEuc)
\end{equation*}
for almost all \((\FlowVar,\pEuc)\in\DFlowXPlus\VField\).

We now recall the notions of distributional derivative
and Lie derivative of a function
\(\funcA\in\LOneLoc\OpenE\) with respect to the vector field \(\VField\).

Following e.g. \cite{article:FSSC:ApproxEmbed:MR1448000}, we say that
\(\funcB\in\LOneLoc{\OpenE}\) is the \emph{distributional derivative}
of \(\funcA\) with respect to \(\VField\)
if for each test function \(\TestFEuc\in C_0^\infty(\OpenE)\) 
\begin{equation*}
	\int_\OpenE\funcB(\pEuc)\TestFEuc(\pEuc)\,d\pEuc
	=-\int_\OpenE\funcA(\pEuc)\VField\TestFEuc(\pEuc)\,d\pEuc
	-\int_\OpenE\funcA(\pEuc)\TestFEuc(\pEuc)\Divergence\VFieldE(\pEuc)\,d\pEuc,
\end{equation*}
where
\begin{equation*}
	\VField\TestFEuc(\pEuc)
	=
	\VFComp_1(\pEuc)\dfrac{\partial\TestFEuc(\pEuc)}{\partial\pEuc_1}
	+\cdots+\VFComp_\Dim(\pEuc)\dfrac{\partial\TestFEuc(\pEuc)}{\partial\pEuc_\Dim}	
\end{equation*}
and
\begin{equation*}
	\Divergence\VFieldE
	\definedby
	\dfrac{\partial\VFComp_1}{\partial\pEuc_1}
	+\cdots+
	\dfrac{\partial\VFComp_\Dim}{\partial\pEuc_\Dim}
\end{equation*}
is the \emph{divergence} of the vector field \(\VField\),
which is a locally bounded functions defined almost everywhere
on \(\OpenE\) by the classical Rademacher theorem.

We say that \(\funcB\in\LOneLoc\OpenE\) is the \emph{Lie derivative}
of \(\funcA\in\LOneLoc\OpenE\) with respect to the vector field \(\VField\)
if
%whenever \([0,\FlowVar_0]\times\SubOpenE\subset\DFlowX\VField\)
%and \(0<\FlowVar<\FlowVar_0\)
%then
%(which may depend on \(\FlowVar	\))
\begin{equation*}
	\funcA\bigl(\VFlow{\VField}{\FlowVar}{\pEuc}\bigr)-\funcA(\pEuc)
%	\DeltaR\VField\FlowVar\funcA(\pEuc)%\funcA\bigl(\VFlow{\VField}{\FlowVar}{\pEuc}\bigr)-\funcA(\pEuc)
	=\int_0^\FlowVar\funcB\bigl(\VFlow{\VField}{\FlowVarB}{\pEuc}\bigr)\, d\FlowVarB,
\end{equation*}
or equivalently
\begin{equation*}
	\DeltaR\VField\FlowVar\funcA(\pEuc)
	=
	\MeanOp{\VFieldE}{\FlowVar}{\funcB}(\pEuc)
\end{equation*}
for almost all \((\FlowVar,\pEuc)\in\DFlowXPlus\VField\).

%\begin{equation*}
%	\funcA\bigl(\VFlow{\VField}{\FlowVar}{\pEuc}\bigr)-\funcA(\pEuc)
%	=\int_0^\FlowVar\funcB\bigl(\VFlow{\VField}{\FlowVarB}{\pEuc}\bigr)\, d\FlowVarB
%\end{equation*}

%We will prove that in this case %for almost all \(\pEuc\in\OpenE\)
%for each relatively compact open set \(\OpenE'\sset\OpenE\)
%\begin{equation*}
%	\lim_{\FlowVar\to0^+}\int_{\OpenE'}
%	\abs{\dfrac{\funcA\bigl(\VFlow{\VField}{\FlowVar}{\pEuc}\bigr)-\funcA(\pEuc)}{\FlowVar}-\funcB(\pEuc)}\,d\pEuc=0,
%\end{equation*}
%that is
%\begin{equation}\label{eq::LieAlmostEverywhere}
%	\lim_{\FlowVar\to0^+}
%	%\dfrac{\funcA\circ\VFlowMap{\VField}{\FlowVar}-\funcA}{\FlowVar}=\funcB
%	\FlowVar^{-1}\bigl(\funcA\circ\VFlowMap{\VField}{\FlowVar}-\funcA\bigr)=\funcB
%\end{equation}
%with respect to the topology of \(\LOneLoc\OpenE\).
%%\begin{equation}\label{eq::LieAlmostEverywhere}
%%	\funcB(\pEuc) = \lim_{\FlowVar\to0^+}\dfrac{\funcA\bigl(\VFlow{\VField}{\FlowVar}{\pEuc}\bigr)-\funcA(\pEuc)}{\FlowVar}.
%%\end{equation}
%%converge when
%It follows that the Lie derivative \(\funcB\), when exists,
%depends uniquely on the function \(\funcA\).

%Clearly the distributional derivative \(\funcB\), if exists,
%it is uniquely determined on the function \(\funcA\).

In this paper, we adopt the notations \(L_\VField\funcA\)
and \(\VField\funcA\) for the Lie derivative
and the distributional derivative respectively
of the function \(\funcA\)
with respect to the vector field \(\VField\).

It is quite easy to prove that if \(\funcA\in C^1(\OpenE)\)
then the function
\begin{equation*}
	\funcB=\VFComp_1\dfrac{\partial\funcA}{\partial\pEuc_1}
	+\cdots+\VFComp_\Dim\dfrac{\partial\funcA}{\partial\pEuc_\Dim}.
\end{equation*}
is continuous on \(\OpenE\) %,satisfies \eqref{eq::LieAlmostEverywhere}
and is both the Lie derivative
and the distributional derivative of the function \(\funcA\).

When \(\funcA\) is not smooth, the relation between \(L_\VField\funcA\)
and \(\VField\funcA\) seems to be not trivial.

For example, in the paper \cite{article:BonfiglioliUguzzoniMR2250624} the authors
prove that if \(\funcA\) and \(L_\VField\funcA\) are both continuous
then \(L_\VField\funcA=\VField\funcA\) under the further hypothesis
that \(\VField\in C^1(\OpenE)\).

Here the first result of this paper.

\begin{theorem}\label{stm::MainThm}
Let \(\VField\) be a locally Lipshitzian vector field on the open set \(\OpenE\subset\RR^\Dim\).
%Then the function \(\funcA\) admits the Lie derivative
%\(L_\VField\funcA\in\LOneLoc\OpenE\)
%if, and only if,
%it admits the distributional derivative \(\VField\funcA\in\LOneLoc\OpenE\)
%and in this case
%\begin{equation*}
%	L_\VField\funcA = \VField\funcA.
%\end{equation*}
%
%and %formula \label{eq::LieAlmostEverywhere}
%\begin{equation}\label{eq::LieAlmostEverywhere}
%	\lim_{\FlowVar\to0^+}\int_{\OpenE'}
%	\abs{\dfrac{\funcA\bigl(\VFlow{\VField}{\FlowVar}{\pEuc}\bigr)-\funcA(\pEuc)}{\FlowVar}-L_\VField\funcA(\pEuc)}\,d\pEuc=0,
%\end{equation}
%for each relatively compact open set \(\OpenE'\sset\OpenE\),
%that is
%\begin{equation*}%\label{eq::LieAlmostEverywhere}
%	\lim_{\FlowVar\to0^+}
%	%\dfrac{\funcA\circ\VFlowMap{\VField}{\FlowVar}-\funcA}{\FlowVar}=\funcB
%	\FlowVar^{-1}\bigl(\funcA\circ\VFlowMap{\VField}{\FlowVar}-\funcA\bigr)=\funcB
%\end{equation*}
%holds
%in \(\LOneLoc\OpenE\). %for almost all \(\pEuc\in\OpenE\).
and let \(\funcA, \funcB\in\LOneLoc\OpenE\) be given.

Then the following conditions are equivalent:
\begin{enumerate}[label=\cbRomanLabel]
\item\label{\stmGFXLoc}\(\funcB=\VFieldE\funcA\), that is \(\funcB\) is
the distributional derivative of \(\funcA\);  
\item\label{\stmGFLLoc}\(\funcB=L_\VFieldE\funcA\), that is \(\funcB\) is
the Lie derivative of \(\funcA\);
\item\label{\stmGFZLoc} %we have
\(%\begin{equation*}%\label{eq::LieAlmostEverywhere}
	\lim_{\FlowVar\to0^+}
	\DeltaR\VField\FlowVar\funcA=\funcB
	%\dfrac{\funcA\circ\VFlowMap{\VField}{\FlowVar}-\funcA}{\FlowVar}=\funcB
	%\FlowVar^{-1}\bigl(\funcA\circ\VFlowMap{\VField}{\FlowVar}-\funcA\bigr)=\funcB
\) %\end{equation*}
with respect to the \(\LOneLoc{\OpenE}\) topology, that is
\begin{equation*}%\label{eq::LieAlmostEverywhere}
	\lim_{\FlowVar\to0^+}\int_{\OpenE'}
	\abs{\dfrac{\funcA\bigl(\VFlow{\VField}{\FlowVar}{\pEuc}\bigr)-\funcA(\pEuc)}{\FlowVar}-\funcB(\pEuc)}\,d\pEuc=0,
\end{equation*}
for each relatively compact open set \(\OpenE'\sset\OpenE\),
\end{enumerate}
%and in this case
%
%Then the function \(\funcA\) admits the Lie derivative
%\(L_\VField\funcA\in\LOneLoc\OpenE\)
%if, and only if,
%it admits the distributional derivative \(\VField\funcA\in\LOneLoc\OpenE\)
%and in this case
%\begin{equation*}
%	L_\VField\funcA = \VField\funcA.
%\end{equation*}
%Moreover %formula \label{eq::LieAlmostEverywhere}
%\begin{equation}\label{eq::LieAlmostEverywhere}
%	\lim_{\FlowVar\to0^+}\int_{\OpenE'}
%	\abs{\dfrac{\funcA\bigl(\VFlow{\VField}{\FlowVar}{\pEuc}\bigr)-\funcA(\pEuc)}{\FlowVar}-\funcB(\pEuc)}\,d\pEuc=0,
%\end{equation}
%for each relatively compact open set \(\OpenE'\sset\OpenE\),
%that is
%\begin{equation*}%\label{eq::LieAlmostEverywhere}
%	\lim_{\FlowVar\to0^+}
%	\dfrac{\funcA\circ\VFlowMap{\VField}{\FlowVar}-\funcA}{\FlowVar}=\funcB
%	\FlowVar^{-1}\bigl(\funcA\circ\VFlowMap{\VField}{\FlowVar}-\funcA\bigr)=\funcB
%\end{equation*}
%holds
%in \(\LOneLoc\OpenE\). %for almost all \(\pEuc\in\OpenE\).
\end{theorem}

The following theorems give a precise relationship
between distributional derivatives and upper gradients.

\begin{theorem}\label{stm::Main::UGradToX}
If the function \(\funcA\in\LOneLoc\OpenE\)
admits an upper gradient \(\UGradA\in\LOneLoc\OpenE\)
along the vector field \(\VFieldE\)
then it also admits the distributional derivative
\(\VFieldE\funcA\in\LOneLoc\OpenE)\)
and 
\begin{equation*}
	\abs{\VFieldE\funcA}\leq\UGradA
\end{equation*}
almost everywhere on \(\OpenE\).
%with respect to \(\VFieldE\).
\end{theorem}

Conversely, we have:

\begin{theorem}\label{stm::Main::XToUGrad}
	Assume that the function \(\funcA\in\LOneLoc\OpenE\)
	admits the distributional derivative
	\(\VFieldE\funcA\in\LOneLoc\OpenE\)
	with respect to \(\VFieldE\).
	Then \(\abs{\VFieldE\funcA}\)
	is an upper gradient of \(\funcA\)
	along the vector field \(\VFieldE\)
	and every upper gradient \(\UGradA\in\LOneLoc\OpenE\)
	of \(\funcA\)
	along the vector field \(\VFieldE\)
	satisfies
	\begin{equation*}
		\abs{\VFieldE\funcA}\leq\UGradA
	\end{equation*}
	almost everywhere on \(\OpenE\).
\end{theorem}

As consequence of the theorems above, we also obtain:

\begin{theorem}\label{stm::Main::UGradToXSystem}
	Let
	\begin{equation*}
		\VFieldE_\VFIMin,\ldots,\VFieldE_\VFIMax
	\end{equation*}
	be locally Lipshitzian vector fields on the open set \(\OpenE\subset\RR^\Dim\).
	%Let \(\funcA, \UGradA\in\LOneLoc{\OpenE}\) be given.
	Assume that the function \(\UGradA\in\LOneLoc{\OpenE}\) is an upper gradient
	of the function \(\funcA\in\LOneLoc{\OpenE}\) along the vector field % of the form
	\begin{equation*}
		\sum_{\VFIndex=\VFIMin}^{\VFIMax}\VConstX_\VFIndex\VFieldE_\VFIndex
	\end{equation*}
	whenever
	%where the constants \(\VConstX_\VFIMin,\ldots,\VConstX_\VFIndex\) stisfies
	\begin{equation*}
		\sum_{\VFIndex=\VFIMin}^{\VFIMax}\VConstX_\VFIndex^2\leq1.
	\end{equation*}
	Then the distributional derivatives
	\( %\begin{equation*}
	\VFieldE_\VFIMin\funcA,\ldots,\VFieldE_\VFIMax\funcA
	\) %\end{equation*}
	are locally integrable and satisfies
	\begin{equation*}
		\sum_{\VFIndex=\VFIMin}^{\VFIMax}(\VFieldE_\VFIndex\funcA)^2\leq\UGradA^2
	\end{equation*}
	almost everywhere on \(\OpenE\).
\end{theorem}

The above theorem extends Theorem 11.7 of
\cite{article:HajlaszKoskela:SobolevMetPoincare},
where the same assertion is proved requiring
the continuity of the function \(\funcA\)
with respect to the Euclidean topology.

}%end unit
%}{}
\subsection{Main ideas}
% !TeX encoding = UTF-8
% !TeX spellcheck = en_GB
{%begin unit
Let \(\VFieldE\) be a complete Lipshitzian vector field on the open set \(\OpenE\subset\RR^\Dim\)
with associated flow \(\VFlow{\VFieldE}{\FlowVar}{\pEuc}\)
\begin{equation*}
	\RR\times\OpenE\ni(\FlowVar,\pEuc)\mapsto\VFlow{\VFieldE}{\FlowVar}{\pEuc}\in\OpenE.
\end{equation*}
%Let denote by \(D\VFlow{\VField}{\FlowVar}{\pEuc}\)
%the Jacobian matrix of the map
%\begin{equation*}
%	\OpenE\ni\pEuc\mapsto\VFlow{\VField}{\FlowVar}{\pEuc}\in\OpenE
%\end{equation*}
%and set
%\begin{equation*}
%	\VFJac\VFieldE\FlowVar\pEuc
%	=\det D\VFlow{\VField}{\FlowVar}{\pEuc}.
%\end{equation*}

The key results of this paper are the following:
\begin{itemize}
\item
for an \emph{arbitrary} function \(\funcA\in \LOneLoc{\OpenE}\)
the difference quotients
\(%\begin{equation*}
	\DeltaR\VField\FlowVar\funcA%(\pEuc)
	%\definedby
	%\dfrac{\funcA\bigl(\VFlow{\VField}{\FlowVar}{\pEuc}\bigr)-\funcA(\pEuc)}{\FlowVar}
\) %\end{equation*}
\emph{always} converge as \(\FlowVar\to0\) (as distribution) to the distribution
\begin{equation*}
	C_0^\infty(\OpenE)\ni\TestFEuc\mapsto
		-\int_\OpenE\funcA\VFieldE\TestFEuc\,d\pEuc
		-\int_\OpenE\funcA\TestFEuc\Divergence\VFieldE\,d\pEuc.
\end{equation*}
(see Proposition \ref{stm::DQLim});

\item if \(\funcA\in L^1(\OpenE)\) and \(\FlowVar\in\RR\) then
\(\funcA\circ\VFlow{\VFieldE}{\FlowVar}{\cdot}\in L^1(\OpenE)\)
and the formula
\begin{equation*}
	\RR\times L^1(\OpenE)\ni(\FlowVar,\funcA)\mapsto
	\GPFlow\VField\FlowVar\funcA\definedby\funcA\circ\VFlowW{\VField}{\FlowVar}\in L^1(\OpenE)
\end{equation*}
defines a joitly continuos (i. e. \(c_0\)) one parameter group
on the Banach space \(L^1(\OpenE)\) (Theorem \ref{stm::FlowToGroup});

\item if \(\GFZ:\Dom(\GFZ)\subset L^1(\OpenE)\to L^1(\OpenE)\)
is the infinitesimal generator of such a semigroup then
\(\funcA\in\Dom(\GFZ)\) if, and only if, \(\VFieldE\funcA\in L^1(\OpenE)\)
and in this case \(\VFieldE\funcA=\GFZ\funcA\)
(Theorem \ref{stm::LOne::ZEqLEqX}). 

\end{itemize}

So, all the machinery of \(c_0-\) one parameter (semi)groups is available.

In this work, we will make frequent use of the
Lipshitzian change of variable in the Lebesgue integral
(see e.g. \cite[Theorem 3.2.5]{book:Federer})

}%end unit
\subsection{Content of the paper}
% !TeX encoding = UTF-8
% !TeX spellcheck = en_GB
{%begin unit
In section \ref{section:Preliminaries}
we review all the prerequisites that we need in this paper.

In section \ref{section:VFEuclidead}
we recall some standard properties of Lipshitzian
vector fields and their associated (local) flows.

Section \ref{section:VFComplete} is the core of the paper:
we prove all the main theorems of the paper
assuming that the vector fields involved
are complete and globally Lipshitzian.

In section \ref{section:MainProof} we give the proof of
the main theorems of the paper by using the results of the previous
section with the help of some standard localization arguments.

}%end unit
\section*{Acknowledgement}
% !TeX encoding = UTF-8
% !TeX spellcheck = en_GB
{%begin unit
I would thank Bruno Franchi and Piotr Haj\l asz for some helpful conversations on this topic.
}%end unit

\section{\label{section:Preliminaries}Preliminaries}
\subsection{Notations}
% !TeX encoding = UTF-8
% !TeX spellcheck = en_GB
{%begin unit

In this paper we will adopt the following notations:
\begin{itemize}
\item \(\OpenE\) is an open subset of \(\RR^\Dim\), the \(\Dim-\)dimensional Euclidean space and \(\overline{\OpenE}\) is
	its closure with respect to the Euclidean topology;
\item \(\EucNorm\pEuc\) is the standard Euclidean norm of \(\pEuc\in\RR^\Dim\);
\item \(\abs{E}\) is the Lebesgue measure of the Borel subset \(E\subset\RR^\Dim\);
\item \(\VFieldE\) is a locally Lipshitzian vector field on \(\OpenE\);
\item \(C_0(\OpenE)\) is the space of the compactly supported
	real continuously functions on \(\OpenE\);
\item \(C_0^k(\OpenE)\), \(k=1,2,\ldots,\infty\) is the space of the compactly supported
	real differentiable functions on \(\OpenE\) of class \(C^k\);
\item \(\LipC{\OpenE}\) is the space of the compactly supported
	real lipshitzian functions on \(\OpenE\);
%\item \(\VFSpace\OpenE\) is the space of all the locally Lipshitzian vector field on the open set \(\OpenE\);
\item \(L^p(\OpenE)\), \(1\leq p\leq+\infty\) is the space of the \(p-\)integrable real functions on \(\OpenE\)
         (with respect to the Lebesgue measure);
\item \(\LOneLoc{\OpenE}\) is the space of the locally integrable real functions on the open set \(\OpenE\)
         (with respect to the Lebesgue measure);
\item \(\BanachS\) is a Banach space with norm \(\BNorm{\cdot}\);
\item \(\BanachS^*\) is the dual of the Banach space \(\BanachS\) and \(\BDuality\cdot\cdot\)
	is the duality map between \(\BanachS\) and \(\BanachS^*\).
%\item \(\OpenE\)

\end{itemize}

}%end unit
%\subsection{Localization}
%\input{src/Localization}
%\subsection{Integration on Banach spaces}
%\input{src/Integration}
\subsection{One parameter semigroups}\label{ssec::OPS}
% !TeX encoding = UTF-8
% !TeX spellcheck = en_GB
{%begin unit
\def\BIndex{k}
\def\BDenseL{S}

In this section we recall some basic fact of the theory
of one parameter semigroups on a Banach space.

All we need is contained in the first chapter of 
\cite{book:davies1980one}.

Let \(\BanachS\) be a Banach space with norm \(\BNorm{\cdot}\).
A continuous map
\begin{equation*}
	[0,+\infty[\times\BanachS\ni(\FlowVar,\BElem)
	\mapsto\BSMap_\FlowVar(\BElem)
	\in\BanachS
\end{equation*}
is a (jointly continuous or \(c_0\)) \emph{one parameter semigroup} on the Banach space \(\BanachS\)
if for each \(\BElem\in\BanachS\)
and each \(\FlowVarA, \FlowVarB\in[0,+\infty[\)

\begin{tabular}{rl}
\Case{1}&\(\BSMap_0(\BElem) = \BElem\);\\
\Case{2}&\(\BSMap_\FlowVarA\bigl(\BSMap_\FlowVarB(\BElem)\bigr)
=\BSMap_{\FlowVarA+\FlowVarB}(\BElem)\).
\end{tabular}

%\Case{1} \(\BSMap_0(\BElem) = \BElem\);
%
%\Case{2} \(\BSMap_\FlowVarA\bigl(\BSMap_\FlowVarB(\BElem)\bigr)
%	=\BSMap_{\FlowVarA+\FlowVarB}(\BElem)\)
%
\begin{proposition}\label{stm::BSCtoJC}
\cite[Proposition 1.18]{book:davies1980one}
A map
\begin{equation*}
	[0,+\infty[\times\BanachS\ni(\FlowVar,\BElem)
	\mapsto\BSMap_\FlowVar(\BElem)
	\in\BanachS
\end{equation*}
which satisfies \Case{1} and \Case{2} above
is a \(c_0\) one parameter semigroup on the Banach space \(\BanachS\)
if, and only if,
for each \(\BElem\in\BanachS\)
\begin{equation*}
	\lim_{\FlowVar\to0}\norm{\BSMap_\FlowVar(\BElem)-\BElem}=0.
\end{equation*}
\end{proposition}

\begin{proposition}
\cite[Proposition 1.23]{book:davies1980one}
A map
\begin{equation*}
	[0,+\infty[\times\BanachS\ni(\FlowVar,\BElem)
	\mapsto\BSMap_\FlowVar(\BElem)
	\in\BanachS
\end{equation*}
which satisfies \Case{1} and \Case{2} above
is a \(c_0\) one parameter semigroup on the Banach space \(\BanachS\)
if, and only if,
for each \(\BElem\in\BanachS\) the map
\begin{equation*}
	[0,+\infty[\ni\FlowVar\mapsto\BSMap_\FlowVar(\BElem)\in\BanachS
\end{equation*}
is weakly continuous at \(\FlowVar=0\),
that is, it is continuous when \(\BanachS\)
is endowed with the weak topology.
\end{proposition}

The (infinitesimal) \emph{generator}
\begin{equation*}
	\BGen:\BDomGen{\BGen}\to\BanachS
\end{equation*}
of the semigroup \(\BSMap_\FlowVar\)
is defined by
\begin{equation*}
	\BGen(\BElem)=\lim_{\FlowVar\to0^+}\FlowVar^{-1}\bigl(\BSMap_\FlowVar(\BElem)-\BElem\bigr)
\end{equation*}
where \(\BDomGen{\BGen}\) is the set of \(\BElem\in\BanachS\)
for which the limit exists.

\begin{proposition}\label{stm::GFIsClosed}
\cite[Lemma 1.1, Lemma 1.5]{book:davies1980one}
The infinitesimal generator \(\BGen\) of the semigroup \(\BSMap_\FlowVar\)
is a closed densely defined linear operator , that is,
the domain \(\BDomGen{\BGen}\) is dense in \(\BanachS\)
and the graph
\begin{equation*}
	\bigl\{(\BElem,\BElemZ)\in\BanachS\times\BanachS\mid\BElem\in\BDomGen{\BGen},\ \BGen(\BElem)=\BElemZ\bigr\}
\end{equation*}
is (a linear subspace) closed in \(\BanachS\times\BanachS\).
\end{proposition}

Let denote by \(\BanachS^*\) be the dual space of \(\BanachS\)
and by \(\BDuality{\cdot}{\cdot}\) the duality map
between \(\BanachS\) and \(\BanachS^*\).

\begin{proposition}\label{stm::GFWStar}
\cite[Theorem 1.24]{book:davies1980one}
Let \(\BElem\) and \(\BElemZ\) be elements of \(\BanachS\).
If there is a subset \(\BDenseL\subset\BanachS^*\) which is \weakstar dense in \(\BanachS^*\)
and a sequence of positive numbers \(\FlowVar_\BIndex\to0\)
such that
\begin{equation*}
	\lim_{\BIndex\to+\infty}\FlowVar_\BIndex^{-1}\bigl\langle\BSMap_{\FlowVar_\BIndex}(\BElem)-\BElem,\BStarElem\bigr\rangle
	=\bigl\langle\BElemZ,\BStarElem\bigr\rangle
\end{equation*}
for each \(\BStarElem\in\BDenseL\)
then \(\BElem\in\BDomGen{\BGen}\) and
\(\BGen(\BElem)=\BElemZ\).

\end{proposition}

\ifdef{\cbFullBase}{%}{}

A vector subspace \(\SCore\subset\BDomGen{\BGen}\) is a \emph{core} for \(\BGen\)
if for all \(\BElem\in\BDomGen{\BGen}\) there exists
a sequence \(\BElem_\BIndex\in\SCore\)
such that \(\BElem_\BIndex\to\BElem\)
and
\(\BGen(\BElem_\BIndex)\to\BGen(\BElem)\)
with respect to the topology induced by the norm of \(\BanachS\).

\begin{proposition}\label{stm::SPCOre}
\cite[Theorem 1.9]{book:davies1980one}
If \(\SCore\subset\BDomGen{\BGen}\) is dense in \(\BanachS\)
and invariant under the semigroup \(\BSMap_\FlowVar\),
that is
\begin{equation*}
	\BSMap_\FlowVar(\SCore)\subset\SCore
\end{equation*}
for each \(\FlowVar\geq0\)
then \(\SCore\) is a core for \(\BGen\).
\end{proposition}

}{}

In this paper we will consider mainly
one parameter semigroups which also are \emph{groups}
that is the operators \(\BSMap_\FlowVar\) are defined for each \(\FlowVar\in\RR\).

The following proposition
\begin{proposition}\label{stm::GPversusSGP}
\cite[Theorem 1.14]{book:davies1980one}
\(\BGen\) is a generator of a \(c_0-\)one parameter group of operators
if, and only if, \(\BGen\) and \(-\BGen\) are both generators of
\(c_0-\)one parameter semigroups of operators.
\end{proposition}
\noindent
often allows us to study a one parameter group
only for \(\FlowVar\geq0\).

}%end unit
\subsection{Banach algebras}\label{ssec::BAlg}
% !TeX encoding = UTF-8
% !TeX spellcheck = en_GB
{%begin unit
\def\BAlg{A}
\def\BElemA{a}
\def\BElemB{b}
\def\UBall{S}
\def\BMap{u}
\def\MatrixA{X}
\def\BELoc{\BElemB}
\def\BEInv{\BElemA}
A \emph{Banach algebra} \(\BAlg\) is a
Banach space which also is a ring where the product
\(\BAlg\times\BAlg\ni(\BElemA,\BElemB)\mapsto\BElemA\BElemB\in\BAlg\)
satisfies
\begin{equation*}
	\norm{\BElemA\BElemB}\leq\norm{\BElemA}\norm{\BElemB}.
\end{equation*}

Assume that the Banach algebra \(\BAlg\) has a unit \(1\in\BAlg\) as a ring
%and \(\norm{1}=1\)
.
If \(\BELoc\in\BAlg\) and \(\norm{\BELoc}<1\)
then \(1+\BELoc\in\BAlg\)
%each element given by \(1+\BELoc\in\BAlg\) where \(\norm{\BELoc}<1\)
is invertible in \(\BAlg\) with inverse given by the ``geometrical series''
\begin{equation*}
	1 - \BELoc + \BELoc^2 - \BELoc^2 + \cdots.
\end{equation*}
In other words, each elment of the ``ball''
\begin{equation*}
	\UBall\definedby\bigl\{\BElemA\in\BAlg\mid\norm{\BElemA-1}<1\bigr\}
\end{equation*}
is invertible in \(\BAlg\to\RR\).

If \(\BMap:\BAlg\to\RR\) is any continuous map which satisfies
\(\BMap(\BEInv)\neq0\) whenever \(\BEInv\)
is invertible in \(\BAlg\) then,
observing that \(\UBall\) is connected,
necessarily either \(\BMap(\UBall)\subset]0,+\infty[\)
or \(\BMap(\UBall)\subset]-\infty,0[\).

In particular, if \(\BMap(1)>0\) then \(\BMap(\UBall)\subset]0,+\infty[\).

Let now
\(\BAlg\) be the algebra
of the real square matrices of order \(\Dim\)
endowed with the operator norm
and let \(\BMap=\det:\BAlg\to\RR\) be the determinant function.

Recalling that \(\BAlg\) is a Banach algebra with unit given by
\(\MatrId\), the identity matrix of order \(\Dim\),
then
the above argument
yields immediately the following result:

\begin{proposition}\label{stm::DetNotZero}
For each real square matrix \(\MatrixA\) of order \(\Dim\)
we have
\begin{equation*}
	\norm{\MatrixA-\MatrId}<1\ \Longrightarrow\ \det\MatrixA>0.
\end{equation*}
\end{proposition}

}%end unit
\subsection{Uniformly integrability}
% !TeX encoding = UTF-8
% !TeX spellcheck = en_GB
{%begin unit
\def\FFunc{f}
\def\FamFunc{\mathcal{F}}
\def\FParam{t}
\def\FIndex{k}
\def\FTest{u}
\def\BSubset{K}
\def\BSubset{E}

Let \(\OpenE\subset\RR^\Dim\) be an open set.
We recall that a family of functions \(\FamFunc\subset L^1(\OpenE)\)
is \emph{uniformly integrable} (or \emph{equi-integrable})
if for each \(\varepsilon>0\) there exists \(\delta>0\)
such that for each \(\FFunc\in\FamFunc\)
\begin{equation*}
	\int_\BSubset\abs{\FFunc(\pEuc)}\,d\pEuc<\varepsilon
\end{equation*}
whenever \(\BSubset\subset\OpenE\) is a Borel subset which satisfies
\begin{equation*}
	\abs{\BSubset}<\delta.
\end{equation*}

%The following proposition is well-known.
%\begin{proposition}
%Let \(\OpenE\in\RR^\Dim\) be an open set
%and let \(\FFunc_\FIndex\in L^1(\OpenE)\)
%be an uniformly integrable sequence of functions
%which converges almost everywhere to a function
%\(%\begin{equation*}
%	\FFunc:\OpenE\to\RR.
%\).%\end{equation*}
%
%Then \(\FFunc\in L^1(\OpenE)\)
%and
%sequence \(\FFunc_\FIndex\) converges to \(\FFunc\)
%in \(L^1(\OpenE)\),
%that is
%\begin{equation*}
%	\lim_{\FIndex\to+\infty}\int_\OpenE\abs{\FFunc_\FIndex(\pEuc)-\FFunc(\pEuc)}\,d\pEuc=0.
%\end{equation*}
%\end{proposition}

We recall the Dunford-Pettis theorem as in \cite{book:FonsecaLeoni:MR2341508}.

\begin{theorem}\label{stm::DunfordPettis}
Let \(\OpenE\in\RR^\Dim\) be an open set
and let \(\FamFunc\subset L^1(\OpenE)\)
a family of functions.

Then the family \(\FamFunc\) is 
that is weakly sequentially precompact in \(L^1(\OpenE)\) if, and only if:
%\begin{enumerate}[label=\((\)\roman*\()\)]
\begin{enumerate}[label=\cbRomanLabel]
\item \(\FamFunc\) is bounded in \(L^1(\OpenE)\);
\item \(\FamFunc\) is uniformly integrable;
\item for each \(\varepsilon>0\) there exists a subset \(\BSubset\subset\OpenE\)
of finite Lebesgue measure such that
for each \(\FFunc\in\FamFunc\)
\begin{equation*}
	\int_{\OpenE\setminus\BSubset}\abs{\FFunc(\pEuc)}\,d\pEuc<\varepsilon.
\end{equation*}
\end{enumerate}
\end{theorem}

%\begin{proposition}
%Let \(\OpenE\in\RR^\Dim\) be an open set
%and let \(\FFunc_\FIndex\in L^1(\OpenE)\)
%be an uniformly integrable sequence of functions
%such that the limit
%\begin{equation*}
%	\lim_{\FIndex\to+\infty}\int_{\OpenE}\FFunc_\FIndex(\pEuc)\FTest(\pEuc)\,d\pEuc
%\end{equation*}
%exists for each \(\FTest\in C^\infty_0(\OpenE)\).
%
%Then the sequence \(\FFunc_\FIndex\) converges to a function \(\FFunc\in  L^1(\OpenE)\)
%with respect to the weak topology of \(L^1(\OpenE)\),
%that is
%\begin{equation*}
%	\lim_{\FIndex\to+\infty}\int_{\OpenE}\FFunc_\FIndex(\pEuc)\FTest(\pEuc)\,d\pEuc
%	=\int_{\OpenE}\FFunc(\pEuc)\FTest(\pEuc)\,d\pEuc
%\end{equation*}
%for each \(\FTest\in L^\infty(\OpenE)\).
%\end{proposition}

}%end unit
%\ifdef{\cbFullBase}{%}{}
%	\subsection{Convolution estimates}
%	\input{src/LipVFConvolution}
%}{}

\section{\label{section:VFEuclidead}Vector Fields on Euclidean Spaces}
% !TeX encoding = UTF-8
% !TeX spellcheck = en_GB
{%begin unit
The results in the first three subsections below
are known but very scattered in the literature.
For the sake of completeness,
we will give all the proofs.

}%end unit
\subsection{Basic estimates}
% !TeX encoding = UTF-8
% !TeX spellcheck = en_GB
{%begin unit

Let \(\VFieldE\) be a vector field on the open set \(\OpenE\subset\RR^\Dim\)
with associated (local) flow \(\VFlow{\VFieldE}{\FlowVar}{\pEuc}\)
\begin{equation*}
	\RR\times\OpenE\supset\DFlowX\VFieldE\ni(\FlowVar,\pEuc)\mapsto\VFlow{\VFieldE}{\FlowVar}{\pEuc}\in\OpenE.
\end{equation*}

We recall that if the vector field is Lipschitzian with Lipshitz constant
\(\LipVF\),  that is
\begin{equation}\label{AppendixVF::Lip}
	\EucNorm{\VFieldE(\pEucA)-\VFieldE(\pEucB)}
	\leq
	\LipVF\EucNorm{\pEucA-\pEucB}
\end{equation}
for each \(\pEucA,\pEucB\in\OpenE\),
where \(\abs{\cdot}\) is the standard Euclidean norm on \(\RR^\Dim\),
then the Gronwall inequality implies that the associated flow
\begin{equation*}
	\RR\times\OpenE\supset\DFlowX\VFieldE\ni(\FlowVar,\pEuc)\mapsto\VFlow{\VField}{\FlowVar}{\pEuc}\in\OpenE.
\end{equation*}
satisfies the estimate
\begin{equation}\label{AppendixVF::Grom}
	\EucNorm{\VFlow{\VField}{\FlowVar}{\pEucA}-\VFlow{\VField}{\FlowVar}{\pEucB}}
	\leq
	e^{\LipVF\abs{\FlowVar}}\EucNorm{\pEucA-\pEucB}.
\end{equation}

%The the Gronwall inequality \eqref{AppendixVF::Grom} implies that
It follows that for each \(\FlowVar\in\RR\) the map
\begin{equation*}
	\pEuc\mapsto\VFlow{\VFieldE}{\FlowVar}{\pEuc}
\end{equation*}
is Lipshitzian and 
we denote by \(\VFJac\VFieldE\FlowVar\pEuc\)
the determinant of its Jacobian matrix,
that is, if
\begin{equation*}
	\VFlow{\VField}{\FlowVar}{\pEuc}=
	\bigl(
	\VFlow{1}{\FlowVar}{\pEuc},
	\ldots,
	\VFlow{\Dim}{\FlowVar}{\pEuc}
	\bigr)
\end{equation*}
then
\begin{equation*}
	\VFJac\VFieldE\FlowVar\pEuc
	=
	\det\left(
	\dfrac{\VFlow{i}{\FlowVar}{\pEuc}}{\partial\pEuc_j}
	\right).
\end{equation*}

\begin{proposition}\label{stm::AdvEstimate}
Let \(\VFieldE\) be a lipschitzian vector field
on an open set \(\OpenE\subset\RR^\Dim\)
with Lispshitz constant \(\LipVF\)
and
let
\begin{equation*}
	\RR\times\OpenE\supset\DFlowX\VFieldE\ni(\FlowVar,\pEuc)\mapsto\VFlow{\VField}{\FlowVar}{\pEuc}\in\OpenE.
\end{equation*}
be the associated flow. %\(\VFlow{\VField}{\FlowVar}{\pEuc}\).

Then the function 
\begin{equation*}
	\DFlowX\VFieldE\ni(\FlowVar,\pEuc)\mapsto
	\VFAdvance(\pEuc,\FlowVar)\definedby\VFlow{\VField}{\FlowVar}{\pEuc}-\pEuc
	\in\RR^\Dim
\end{equation*}
satisfies the estimate
\begin{equation}\label{eq::AdvEstimate}
	\abs{\VFAdvance(\pEucA,\FlowVar)-\VFAdvance(\pEucB,\FlowVar)}
	\leq
	%\LipVF\abs{\FlowVar}e^{\LipVF\abs{\FlowVar}}\abs{\pEucA-\pEucB}.
	(e^{\LipVF\abs{\FlowVar}}-1)\EucNorm{\pEucA-\pEucB}.
\end{equation}
\end{proposition}

\begin{proof}{
\def\varA{\FlowVar}
\def\varB{xxx}
\def\varI{s}
\def\auxF{g}

Let \(\varA\in\RR\) and \(\pEucA, \pEucB\in\OpenE\) such that
\((\varA,\pEucA), (\varA,\pEucB)\in\DFlow\).

If \(\varA=0\) then the estimate is obviously satisfied.

Assume that \(\varA>0\).
Then
\begin{equation*}
	\VFlow{\VField}{\varA}{\pEucA}
	=\pEucA+\int_0^\varA\VFieldE\bigl(\VFlow{\VField}{\varI}{\pEucA}\bigr)\,d\varI
\end{equation*}
and
\begin{equation*}
	\VFlow{\VField}{\varA}{\pEucB}
	=\pEucB+\int_0^\varA\VFieldE\bigl(\VFlow{\VField}{\varI}{\pEucB}\bigr)\,d\varI.
\end{equation*}

It follows that
\begin{equation*}
	\VFAdvance(\pEucA,\FlowVar)-\VFAdvance(\pEucB,\FlowVar)
	=
	\int_0^\varA
		\Bigl(\VFieldE\bigl(\VFlow{\VField}{\varI}{\pEucA}\bigr)
		-\VFieldE\bigl(\VFlow{\VField}{\varI}{\pEucB}\bigr)
		\Bigr)\,d\varI.
\end{equation*}
Since \(0\leq\varI\leq\varA\), then the inequalities \eqref{AppendixVF::Lip} and \eqref{AppendixVF::Grom}
imply that
\begin{equation*}
	\abs{\VFieldE\bigl(\VFlow{\VField}{\varI}{\pEucA}\bigr)
	-\VFieldE\bigl(\VFlow{\VField}{\varI}{\pEucB}\bigr)}
	\leq\LipVF e^{\LipVF\varI}\abs{\pEucA-\pEucB}
	%\leq\LipVF e^{\LipVF\varA}\abs{\pEucA-\pEucB}
\end{equation*}
and hence
\begin{equation*}
	\EucNorm{\VFAdvance(\pEucA,\FlowVar)-\VFAdvance(\pEucB,\FlowVar)}
	\leq\int_{0}^{\FlowVar}\LipVF e^{\LipVF\varI}\EucNorm{\pEucA-\pEucB}\,d\varI
	=(e^{\LipVF\FlowVar}-1)\EucNorm{\pEucA-\pEucB},
\end{equation*}
%\begin{equation*}
%	\abs{\VFAdvance(\pEucA,\FlowVar)-\VFAdvance(\pEucB,\FlowVar)}
%	\leq
%	\LipVF\FlowVar e^{\LipVF\FlowVar}\abs{\pEucA-\pEucB},
%\end{equation*}
which is the desired estimate when \(\FlowVar>0\).

The proof when \(\FlowVar<0\)
%follows considering the flow associated to the vector field \(-\VFieldE\).
is similar.

}\end{proof}

%Fix \(\varA\in\RR\), \(\pEuc\in\RR^\Dim\).
%The function
%\begin{equation*}
%	[0,1]\ni\varB\mapsto\auxF(\varB)\definedby\VFlow{\VField}{\varB\varA}{\pEuc}
%\end{equation*}
%is of class \(C^1\) on the interval \([0,1]\) and satisfies
%\begin{eqnarray*}
%	&&\auxF(0)=\pEuc,\\
%	&&\auxF(1)=\VFlow{\VField}{\varA}{\pEuc},\\
%	&&\auxF'(\varB)=\varA\VField\bigl(\VFlow{\VField}{\varB\varA}{\pEuc}\bigr)
%\end{eqnarray*}
%and hence
%\begin{eqnarray*}
%	\VFAdvance(\pEuc,\varA)
%	=\auxF(1)-\auxF(0)
%	=\int_{0}^{1}\auxF'(\varB)\,d\varB
%	=\varA\int_{0}^{1}\VField\bigl(\VFlow{\VField}{\varB\varA}{\pEuc}\bigr)\,d\varB
%\end{eqnarray*}

}%end unnit
\subsection{Jacobian estimates}
% !TeX encoding = UTF-8
% !TeX spellcheck = en_GB
{%begin unit
\def\KStep{k}

Let \(\VField\)
be a complete Lipshitzian vector field
on an open set \(\OpenE\subset\RR^\Dim\)
with Lipshitz constant \(\LipVF\)
and associated flow
\begin{equation*}
	\RR\times\OpenE\ni(\FlowVar,\pEuc)\mapsto\VFlow{\VField}{\FlowVar}{\pEuc}\in\OpenE.
\end{equation*}

We recall that \(D\VFlow{\VField}{\FlowVar}{\pEuc}\)
denotes the Jacobian matrix of the map
\begin{equation*}
	\OpenE\ni\pEuc\mapsto\VFlow{\VField}{\FlowVar}{\pEuc}\in\OpenE
\end{equation*}
and
\begin{equation*}
	\VFJac\VFieldE\FlowVar\pEuc
	=\det D\VFlow{\VField}{\FlowVar}{\pEuc}.
\end{equation*}

The purpose of this section is to prove that
for each \(\FlowVar\in\RR\)
%the Jacobian determinant \(\VFlow{\VField}{\FlowVar}{\cdot}\)
%satisfies
the estimate
\begin{equation}\label{stm::VFJExpEstimate}
	e^{-\Dim\LipVF\abs{\FlowVar}}
	\leq
	\VFJac\VFieldE\FlowVar\pEuc
	\leq
	e^{\Dim\LipVF\abs{\FlowVar}}
\end{equation}
do hold for almost all \(\pEuc\in\OpenE\).

Let denote by  \(\MatrNorm{\cdot}\)
the the operator norm on the space of
the real square matrix of order \(\Dim\).

We begin observing that the estimate \eqref{AppendixVF::Grom}
implies that if \(\VFlow{\VField}{\FlowVar}{\cdot}\)
is differentiable at \(\pEuc\in\OpenE\) then
\begin{equation*}
	\MatrNorm{D\VFlow{\VField}{\FlowVar}{\pEuc}}
	\leq e^{\LipVF\abs{\FlowVar}}
\end{equation*}
and hence
\begin{equation*}
	\abs{\VFJac\VFieldE\FlowVar\pEuc}
	\leq
	\MatrNorm{D\VFlow{\VField}{\FlowVar}{\pEuc}}^\Dim
	\leq e^{\Dim\LipVF\abs{\FlowVar}}.
\end{equation*}

A standard argument implies that
if \(\VFlow{\VField}{\FlowVar}{\cdot}\)
is differentiable at \(\pEuc\in\OpenE\)
then \(\VFlow{\VField}{-\FlowVar}{\cdot}\)
is differentiable at \(\VFlow{\VField}{\FlowVar}{\pEuc}\)
and the matrices 
\(D\VFlow{\VField}{\FlowVar}{\pEuc}\)
and
\(D\VFlowB{\VField}{-\FlowVar}{\VFlow{\VField}{\FlowVar}{\pEuc}}\)
are one the inverse of the other.

Then we have equality
\begin{equation*}
	\VFJac\VFieldE\FlowVar\pEuc
	\VFJacB\VFieldE{-\FlowVar}{\VFlow{\VField}{\FlowVar}{\pEuc}}=1
\end{equation*}
and hence
\begin{equation*}
	\abs{\VFJac\VFieldE\FlowVar\pEuc}
	=
	\abs{\VFJacB\VFieldE{-\FlowVar}{\VFlow{\VField}{\FlowVar}{\pEuc}}}^{-1}
	\geq
	e^{-\Dim\LipVF\abs{\FlowVar}}
	.
\end{equation*}

%First we prove that
In order to prove the estimates \eqref{stm::VFJExpEstimate}
it then suffices to show that
for almost all \(\pEuc\in\OpenE\)
\begin{equation*}
	\VFJac\VFieldE\FlowVar\pEuc>0.
\end{equation*}
It follows from the estimate \eqref{eq::AdvEstimate}
of Proposition \ref{stm::AdvEstimate}
that 
\begin{equation*}
	\MatrNorm{D\VFlow{\VField}{\FlowVar}{\pEuc}-\MatrId}
	\leq e^{\LipVF\abs{\FlowVar}}-1
\end{equation*}
whenever the map \(\VFlow{\VField}{\FlowVar}{\cdot}\) 
is differentiable at \(\pEuc\).
%for almost all \(\pEuc\in\OpenE\).
Here, as usual, \(\MatrId\) denotes the identity matrix of order \(\Dim\).

Observe now that if
\begin{equation*}
	\abs{\FlowVar}<\dfrac{1}{\LipVF}\log2
\end{equation*}
one has
\begin{equation*}
	\MatrNorm{D\VFlow{\VField}{\FlowVar}{\pEuc}-\MatrId}
	\leq e^{\LipVF\abs{\FlowVar}}-1<1
\end{equation*}
%and Lemma xxx yields
%from which it follows that
and hence Proposition \ref{stm::DetNotZero} yields
\begin{equation*}
	\VFJac\VFieldE\FlowVar\pEuc=\det D\VFlow{\VField}{\FlowVar}{\pEuc}>0.
\end{equation*}

Let now \(\FlowVar\in\RR\) be arbitrary.
If \(\KStep\in\NN\) \(\KStep>0\)
for almost all \(\pEuc\in\OpenE\) we have the matrix identity
\begin{equation*}
	D\VFlow{\VField}{\FlowVar}{\pEuc}=
	D\VFlowB{\VField}{\frac{\FlowVar}{\KStep}}{\VFlow{\VField}{\frac{\KStep-1}{\KStep}\FlowVar}{\pEuc}}
		\cdots
	D\VFlowB{\VField}{\frac{\FlowVar}{\KStep}}{\VFlow{\VField}{\frac{\FlowVar}{\KStep}}{\pEuc}}
	D\VFlow{\VField}{\frac{\FlowVar}{\KStep}}{\pEuc}
\end{equation*}
which implies that
\begin{equation*}
	\VFJac\VFieldE\FlowVar\pEuc=
	\prod_{j=0}^{\KStep-1}\VFJacB\VFieldE{\frac{\FlowVar}{\KStep}}{\VFlow{\VField}{\frac{j}{\KStep}}{\pEuc}}.
\end{equation*}
If we choose \(\KStep>0\)
big enough in such a way that
\begin{equation*}
	\dfrac{\abs{\FlowVar}}{\KStep}<\dfrac{1}{\LipVF}\log2.
\end{equation*}
then we already know that
for each \(j=0,\ldots,\KStep-1\)
\begin{equation*}
	\VFJacB\VFieldE{\frac{\FlowVar}{\KStep}}{\VFlow{\VField}{\frac{j}{\KStep}}{\pEuc}}>0
\end{equation*}
and hence we also have
\begin{equation*}
	\VFJac\VFieldE\FlowVar\pEuc=
\prod_{j=0}^{\KStep-1}\VFJacB\VFieldE{\frac{\FlowVar}{\KStep}}{\VFlow{\VField}{\frac{j}{\KStep}}{\pEuc}}>0,
\end{equation*}
as desired.

\begin{remark}
For some better estimates 
related to \eqref{stm::VFJExpEstimate}
%(without proofs)
see, e. g., the introduction of \cite{article:DiPernaLions89:MR1022305}.

\end{remark}

}%end unit
%\subsection{Determinant estimates}
%\input{src/AppendixVFDet}
\subsection{The divergence}
% !TeX encoding = UTF-8
% !TeX spellcheck = en_GB

{%begin unit
\def\JDiff{\dfrac{\VFJac\VFieldE\FlowVar\pEuc-1}{\FlowVar}}

Let
\begin{equation*}
	\VFieldE=\VFComp_1\dfrac{\partial}{\partial\pEuc_1}
	+\cdots+\VFComp_\Dim\dfrac{\partial}{\partial\pEuc_\Dim}
\end{equation*}
be a lipshitzian vector field
on an open set \(\OpenE\subset\RR^\Dim\)
with associated flow
\begin{equation*}
	\RR\times\OpenE\supset\DFlow\ni(\FlowVar,\pEuc)\mapsto\VFlow{\VField}{\FlowVar}{\pEuc}\in\OpenE.
\end{equation*}
Then the components \(\VFComp_1,\ldots,\VFComp_\Dim\)
are lipshitzian functions on \(\OpenE\)
and hence
they are almost everywhere differentiable (with respect to the Lebesgue measure)
by the Rademacher theorem.

The \emph{divergence} of the vector field \(\VFieldE\)
\begin{equation*}
	\Divergence\VFieldE
	\definedby
	\dfrac{\partial\VFComp_1}{\partial\pEuc_1}
	+\ldots+
	\dfrac{\partial\VFComp_\Dim}{\partial\pEuc_\Dim}
\end{equation*}
is then a bounded Borel function defined almost everywhere
on \(\OpenE\).

We recall that \(\VFJac\VFieldE\FlowVar\pEuc\) denotes
the determinant of the Jacobian matrix of the map
\begin{equation*}
	\pEuc\mapsto\VFlow{\VFieldE}{\FlowVar}{\pEuc}.
\end{equation*}

\begin{proposition}\label{stm::VFJacToDiv}
Let \(\VField\)
be a complete Lipshitzian vector field
on an open set \(\OpenE\subset\RR^\Dim\)
with associated flow
\begin{equation*}
	\RR\times\OpenE\ni(\FlowVar,\pEuc)\mapsto\VFlow{\VField}{\FlowVar}{\pEuc}\in\OpenE.
\end{equation*}
Then:

\Case{\stmJBouned}
for each \(\FlowVar_0>0\) the functions
\begin{equation*}
	\JDiff,\ -\FlowVar_0\leq\FlowVar\leq\FlowVar_0,\ \FlowVar\neq0
\end{equation*}
are uniformly bounded in \(L^\infty(\OpenE)\);

\Case{\stmJwStarToDiv}
if \(\TestFEuc\in L^1(\OpenE)\) then
\begin{equation}\label{eq::DivLimInt}
	\lim_{\FlowVar\to0}\int_\OpenE\TestFEuc(\pEuc)\JDiff\,d\pEuc
	=\int_\OpenE\TestFEuc(\pEuc)\Divergence\VFieldE(\pEuc)\,d\pEuc,
\end{equation}
that is
\begin{equation*}%\label{eq::DivLimJ}
	\lim_{\FlowVar\to0}\dfrac{\VFJac\VFieldE\FlowVar\cdot-1}{\FlowVar}=\Divergence\VFieldE
\end{equation*}
with respect to the \weakstar topology of \(L^\infty(\OpenE)\).

%\Case{3}
%the limit \eqref{eq::DivLimJ} holds in \(\LOneLoc{\OpenE}\),
%that is,
%\begin{equation*}
%	\lim_{\FlowVar\to0}\int_{\OpenE'}\abs{\JDiff
%		-\Divergence\VFieldE(\pEuc)}d\pEuc
%	=0
%\end{equation*}
%for each relatively compact open subset
%\(\OpenE'\Subset\OpenE\).

\end{proposition}

{
\def\Constant{K}
\begin{proof}
Let \(\FlowVar_0>0\) be fixed
and assume that
\begin{equation*}
	\abs{\FlowVar}\leq\FlowVar_0.
\end{equation*}
It is an elementary consequence of
\eqref{stm::VFJExpEstimate} that
\begin{equation*}
	\abs{\VFJac\VFieldE\FlowVar\pEuc-1}\leq e^{\LipVF\abs{\FlowVar}}-1
\end{equation*}
and hence
\begin{equation*}
	\abs{\JDiff}
	\leq\dfrac{e^{\LipVF\abs{\FlowVar}}-1}{\abs{\FlowVar}}
	\leq\dfrac{e^{\LipVF\FlowVar_0}-1}{\FlowVar_0}
\end{equation*}
almost everywhere on \(\OpenE\).
%
%Let denote by \(D\VFlow{\VField}{\FlowVar}{\pEuc}\)
%the Jacobian matrix of the map
%\begin{equation*}
%	\OpenE\ni\pEuc\mapsto\VFlow{\VField}{\FlowVar}{\pEuc}.
%\end{equation*}
%Then
%\begin{equation*}
%	\VFJac\VFieldE\FlowVar\pEuc
%	=\det D\VFlow{\VField}{\FlowVar}{\pEuc}.
%\end{equation*}
%It follows from the estimate \eqref{eq::AdvEstimate}
%of Proposition \ref{stm::AdvEstimate}
%that if we endow the space of the
%real square matrix of the operator norm \(\MatrNorm{\cdot}\)
%then
%\begin{equation*}
%	\MatrNorm{D\VFlow{\VField}{\FlowVar}{\pEuc}-\MatrId}
%	\leq e^{\LipVF\abs{\FlowVar}}-1
%	\leq e^{\LipVF\abs{\FlowVar_0}}-1
%\end{equation*}
%for almost all \(\pEuc\in\OpenE\).
%Here, as usual, \(\MatrId\) denotes the identity matrix of order \(\Dim\).
%
%Proposition \ref{stm::DetEstimate} then implies that
%\begin{eqnarray*}
%	\abs{\VFJac\VFieldE\FlowVar\pEuc-1}
%	&=&
%	%\det\bigl(\MatrId+(D\VFlow{\VField}{\FlowVar}{\pEuc}-\MatrId)\bigr)
%	\abs{\det\bigl(D\VFlow{\VField}{\FlowVar}{\pEuc}\bigr)-1}
%	\\
%	&\leq&
%	\Constant\MatrNorm{D\VFlow{\VField}{\FlowVar}{\pEuc}-\MatrId}
%	\leq
%	\Constant(e^{\LipVF\abs{\FlowVar}}-1)
%\end{eqnarray*}
%for a suitable constant \(\Constant\)
%and hence when \(\FlowVar\neq0\)
%\begin{equation*}
%	\abs{\JDiff}
%	\leq\Constant\dfrac{e^{\LipVF\abs{\FlowVar}}-1}{\abs{\FlowVar}}
%	\leq\Constant\dfrac{e^{\LipVF\FlowVar_0}-1}{\FlowVar_0}.
%\end{equation*}

This proves (\rmnum{1}),
which together the density 
of \(C_0^1(\OpenE)\)
in \(L^1(\OpenE)\),
%the compactly supported 
implies that 
%It remains to prove \eqref{eq::DivLimJ}.
in order to prove (\rmnum{2})
it suffices
to check \eqref{eq::DivLimInt} when \(\TestFEuc\in C_0^1(\OpenE)\).

So let \(\TestFEuc\in C_0^1(\OpenE)\) and consider the integral
\begin{equation*}
	\int_\OpenE\VFieldE\COneFTest(\pEuc)\,d\pEuc.
\end{equation*}

%Lemma \ref{stm::EucIntPart} implies that
Of course we have
\begin{equation}\label{eq::local::A}
	\int_\OpenE\VFieldE\COneFTest(\pEuc)\,d\pEuc
	=-\int_\OpenE\COneFTest(\pEuc)\Divergence\VFieldE(\pEuc)\,d\pEuc,
\end{equation}
but we also have
\begin{eqnarray*}
	\int_\OpenE\VFieldE\COneFTest(\pEuc)\,d\pEuc
	&=&\int_\OpenE\lim_{\FlowVar\to0}
		\dfrac{\COneFTest\bigl(\VFlow{\VField}{\FlowVar}{\pEuc}\bigr)
			-\COneFTest(\pEuc)}{\FlowVar}
		\,d\pEuc,\\
	&=&\lim_{\FlowVar\to0}\int_\OpenE
		\dfrac{\COneFTest\bigl(\VFlow{\VField}{\FlowVar}{\pEuc}\bigr)
			-\COneFTest(\pEuc)}{\FlowVar}
		\,d\pEuc
\end{eqnarray*}
by the Lebesgue theorem on the dominated convergence.

The change of variable %formula %for Lipshitzian function yields
\begin{equation*}
	\pEuc\mapsto\VFlow{\VField}{-\FlowVar}{\pEuc}
\end{equation*}
gives
\begin{equation*}
	\int_\OpenE\COneFTest\bigl(\VFlow{\VField}{\FlowVar}{\pEuc}\bigr)\,d\pEuc
	=\int_\OpenE\COneFTest(\pEuc)\VFJac\VFieldE{-\FlowVar}\pEuc\,d\pEuc
\end{equation*}
and hence
\begin{eqnarray*}
	\int_\OpenE
	\dfrac{\COneFTest\bigl(\VFlow{\VField}{\FlowVar}{\pEuc}\bigr)
		-\COneFTest(\pEuc)}{\FlowVar}
	\,d\pEuc
	&=&\dfrac{1}{\FlowVar}\left[
		\int_\OpenE\COneFTest\bigl(\VFlow{\VField}{\FlowVar}{\pEuc}\bigr)\,d\pEuc
		-\int_\OpenE\COneFTest(\pEuc)\,d\pEuc
		\right]\\	
	&=&\dfrac{1}{\FlowVar}\left[
		 \int_\OpenE\COneFTest(\pEuc)\VFJac\VFieldE{-\FlowVar}\pEuc\,d\pEuc
		-\int_\OpenE\COneFTest(\pEuc)\,d\pEuc
	\right]\\
	&=&-\int_\OpenE\COneFTest(\pEuc)
			\dfrac{\VFJac\VFieldE{-\FlowVar}\pEuc-1}{-\FlowVar}\,d\pEuc.
\end{eqnarray*}
It follows that
\begin{equation*}%\label{eq::local::B}
	\int_\OpenE\VFieldE\COneFTest(\pEuc)\,d\pEuc
	=\lim_{\FlowVar\to0}\int_\OpenE
	\dfrac{\COneFTest\bigl(\VFlow{\VField}{\FlowVar}{\pEuc}\bigr)
		-\COneFTest(\pEuc)}{\FlowVar}
		\,d\pEuc
	=-\lim_{\FlowVar\to0}\int_\OpenE\COneFTest(\pEuc)
		\dfrac{\VFJac\VFieldE{\FlowVar}\pEuc-1}{\FlowVar}
		\,d\pEuc
\end{equation*}
which compared with \eqref{eq::local::A} %and \eqref{eq::local::B}
yelds the desired equality \eqref{eq::DivLimInt}.

%Assertion \Case{3} is a standard consequence of \Case{\stmJwStarToDiv}.

\end{proof}
}%end proof

\begin{remark}
A proof of \Case{1} of the previous proposition 
is given inside the proof of
\cite[Proposition 2.9]{article:FSSC:ApproxEmbed:MR1448000}.

\end{remark}

}%end unnit
%\section{\label{section:meanOp}The mean operator}
\subsection{Lebesgue points}
% !TeX encoding = UTF-8
% !TeX spellcheck = en_GB
{%begin unit
\def\varI{s}
\def\varB{r}

\def\LPointF{f}
In this subsection we prove 
an auxiliary result on locally integrable functions
that we need in this paper.

\begin{proposition}\label{stm::LipVFLebesgueP}
Let \(\VField\)
be a locally Lipshitzian vector field
on an open set \(\OpenE\subset\RR^\Dim\).

If \(\LPointF\in\LOneLoc{\OpenE}\) then
\begin{equation*}\label{stm::VFLebesgueX}
	\LPointF(\pEuc)=
		\lim_{\FlowVar\to0}
		\MeanOp{\VField}{\FlowVar}{\funcA}(\pEuc)
		=\lim_{\FlowVar\to0}
		\dfrac{1}{\FlowVar}
		\int_0^\FlowVar\LPointF\bigl(\VFlow{\VField}{\varI}{\pEuc}\bigr)\,d\varI
\end{equation*}
for almost all \(\pEuc\in\OpenE\).

\end{proposition}

\begin{proof} {
\def\MLebesgue#1{\abs{#1}}
\def\AuxF{F}
\def\LPoints{E}
\def\LPointsT#1{\LPoints_{#1}}
\def\LPointsP#1{\LPoints^{#1}}
By a standard localization argument
%(see, e. g., Proposition \ref{stm::LocalizationPrinciple})
we may suppose that
\(\LPointF\in L^1(\OpenE))\) and that
the vector field has compact support in \(\OpenE\)
and hence is Lipshitzian and complete.

Consider the the function
\begin{equation*}
	[0,1]\times\OpenE\ni(\FlowVar,\pEuc)\mapsto
	\AuxF(\FlowVar,\pEuc)\definedby
	\LPointF\bigl(\VFlow{\VField}{\FlowVar}{\pEuc}\bigr)
\end{equation*}
and let \(\LPoints\)
the subset of all pairs \((\FlowVar,\pEuc)\in[0,1]\times\OpenE\)
such that \(\FlowVar\)
%for almost all \(\pEuc\in\OpenE\)
%almost all \(\varB\in[0,1]\) is
is a Lebesgue point for \(\AuxF(\cdot,\pEuc)\),
that is
\begin{equation*}
	\AuxF(\FlowVar,\pEuc)=
	\lim_{\varB\to0}
	\dfrac{1}{\varB}
	\int_\FlowVar^{\FlowVar+\varB}\AuxF(\varI,\pEuc)\,d\varI=
	\lim_{\varB\to0}
	\dfrac{1}{\varB}
	\int_\FlowVar^{\FlowVar+\varB}
		\LPointF\bigl(\VFlow{\VField}{\varI}{\pEuc}\bigr)\,d\varI.
\end{equation*}

Using the continuity of the function
\begin{equation*}
	\varB\mapsto\int_\FlowVar^{\FlowVar+\varB}\AuxF(\varI,\pEuc)\,d\varI
\end{equation*}
%with respect to the variable \(\varB\)
it is not difficult to see that \(\LPoints\) is
a Borel subset of \([0,1]\times\OpenE\).

For \(\pEuc\in\OpenE\) and \(\FlowVar\in[0,1]\)
set respectively
\begin{equation*}
\LPointsP{\pEuc}=\bigl\{\FlowVar\in[0,1]\mid(\FlowVar,\pEuc)\in\LPoints\bigr\}
\end{equation*}
and
\begin{equation*}
	\LPointsT{\FlowVar}=\bigl\{\pEuc\in\OpenE\mid(\FlowVar,\pEuc)\in\LPoints\bigr\}.
\end{equation*}

We need to prove that %almost \(\pEuc\in\OpenE\) belongs to 
\(\MLebesgue{\OpenE\setminus\LPointsT{0}}=0\),
where \(\MLebesgue\cdot\) denotes the Lebesgue measure
(both in \(\RR\) and in \(\RR^\Dim\)).

The change of variables
\begin{equation*}
	\pEuc\mapsto\VFlow{\VField}{-\FlowVar}{\pEuc}
\end{equation*}
and (\rmnum{\stmJBouned}) of Proposition \ref{stm::VFJacToDiv}
yield
\begin{equation*}
\begin{split}
	\int_0^1d\FlowVar\int_\OpenE\abs{\AuxF(\FlowVar,\pEuc)}\,d\pEuc
	&=
	\int_0^1d\FlowVar\int_\OpenE\abs{\LPointF\bigl(\VFlow{\VField}{\FlowVar}{\pEuc}\bigr)}\,d\pEuc
	\\
	&=\int_0^1d\FlowVar\int_\OpenE\abs{\LPointF(\pEuc)}\VFJac\VFieldE{-\FlowVar}\pEuc\,d\pEuc
	<\infty
\end{split}
\end{equation*}
being \(\LPointF\in L^1(\OpenE)\) and \(\VFJac\VFieldE{-\FlowVar}\cdot\)
uniformely bounded
in \(L^\infty(\OpenE)\)
with respect to \(\FlowVar\in[0,1]\)
.

The Fubini-Tonelli theorems imply that
\(\AuxF(\FlowVar,\pEuc)\in L^1([0,1]\times\OpenE)\)
and for almost all \(\pEuc\in\OpenE\)
\begin{equation*}
	\AuxF(\cdot,\pEuc)\in L^1([0,1]).
\end{equation*}

Since almost all point of an integrable function on the interval \([0,1]\)
are Lebesgue points we have thus proved that
\begin{equation*}
	\MLebesgue{[0,1]\setminus\LPointsP{\pEuc}}=0
\end{equation*}
for almost all \(\pEuc\in\OpenE\).

%It follow that if we denote by \(\LPoints\)
%the subset of all pairs \((\FlowVar,\pEuc)\in[0,1]\times\OpenE\)
%such that \(\FlowVar\)
%%for almost all \(\pEuc\in\OpenE\)
%%almost all \(\varB\in[0,1]\) is
%is a Lebesgue point for \(\AuxF(\cdot,\pEuc)\),
%that is
%\begin{equation*}
%	\AuxF(\varB,\pEuc)=
%	\lim_{\FlowVar\to0}
%	\dfrac{1}{\FlowVar}
%	\int_\varB^{\varB+\varI}\AuxF(\varI,\pEuc)\,d\varI,
%\end{equation*}
%for almost all \(\pEuc\in\OpenE\)
%the pair \((\FlowVar,\pEuc)\) belongs to \(\LPoints\)
%for almost all \(\varB\in[0,1]\).

Again by the Fubini theorem it follows that
\begin{equation*}
	\MLebesgue{\OpenE\setminus\LPointsT{\FlowVar}}=0
\end{equation*}
for almost all \(\FlowVar\in[0,1]\).

Let fix \(\FlowVar\in[0,1]\) which satisfies
\(\MLebesgue{\OpenE\setminus\LPointsT{\FlowVar}}=0\).
Since
\begin{equation*}
	\AuxF(\FlowVar,\pEuc)
	=\AuxF\bigl(0,\VFlow{\VField}{\FlowVar}{\pEuc}\bigr)
\end{equation*}
and
\begin{equation*}
\begin{split}
	\lim_{\varB\to0}
		\dfrac{1}{\varB}
		\int_\FlowVar^{\FlowVar+\varB}
		\AuxF(\varI,\pEuc)\,d\varI
	&=
	\lim_{\varB\to0}
		\dfrac{1}{\varB}
		\int_\FlowVar^{\FlowVar+\varB}
		\LPointF\bigl(\VFlow{\VField}{\varI}{\pEuc}\bigr)\,d\varI
	\\
	&=
	\lim_{\varB\to0}
		\dfrac{1}{\varB}
		\int_0^{\varB}
		\LPointF\bigl(\VFlow{\VField}{\FlowVar+\varI}{\pEuc}\bigr)\,d\varI
	\\
	&=
		\lim_{\varB\to0}
		\dfrac{1}{\varB}
		\int_0^{\varB}
		\LPointF\bigl(
			\VFlow{\VField}{\varI}{\VFlow{\VField}{\FlowVar}{\pEuc}}
		\bigr)\,d\varI
	\\
	&=
		\lim_{\varB\to0}
		\dfrac{1}{\varB}
		\int_0^{\varB}
		\AuxF\bigl(\varI,\VFlow{\VField}{\FlowVar}{\pEuc}\bigr)\,d\varI
\end{split}
\end{equation*}
it follows that
\begin{equation*}
	(\FlowVar,\pEuc)\in\LPoints
	\ \Longleftrightarrow\ (0,\VFlow{\VField}{\FlowVar}{\pEuc})\in\LPoints.
\end{equation*}
%\begin{equation*}
%	\VFlow{\VField}{\FlowVar}{\LPointsT{\FlowVar}}
%	=\LPointsT{0}.
%\end{equation*}

Since the map
\begin{equation*}
	\OpenE\ni\pEuc\mapsto\VFlow{\VField}{\FlowVar}{\pEuc}\in\OpenE
\end{equation*}
is a bi-Lipshitz homeomorphism it follows that
\begin{equation*}
	\MLebesgue{\OpenE\setminus\LPointsT{0}}=
	{\MLebesgue{\VFlow{\VField}{\FlowVar}{\OpenE\setminus\LPointsT{\FlowVar}}}}=0,
\end{equation*}
as required.

}
\end{proof}

%As immediace corollary we obtain:
A standard argument using the Fubini-Tonelli theorems gives the
following Proposition:

\begin{proposition}\label{stm::LIpVF::PointwiseLieDer}
Let \(\VFieldE\)
be a locally Lipshitzian vector field
on an open set \(\OpenE\subset\RR^\Dim\)
and let
\(\funcA,\funcB\in\LOneLoc{\OpenE}\).

If \(\funcB=L_\VFieldE\funcA\) then
for almost all \(\pEuc\in\OpenE\).
it is possible to find a subset \(\NullSet\subset\RR\),
\(\NullSet=\NullSet(\pEuc)\),
having Lebesgue measure zero such that
\begin{equation*}
	\funcB(\pEuc)=\lim_{\RR\setminus\NullSet\ni\FlowVar\to0}
	\MeanOp{\VField}{\FlowVar}{\funcA}(\pEuc)
	=\lim_{\RR\setminus\NullSet\ni\FlowVar\to0}
	\dfrac{\funcA\bigl(\VFlow{\VField}{\FlowVar}{\pEuc}\bigr)-\funcA(\pEuc)}{\FlowVar}.
\end{equation*}
\end{proposition}

An other usefull consequence of proposition \ref{stm::LipVFLebesgueP} is
the following:

\begin{proposition}\label{stm::AlmostLieToLie}
Let \(\VFieldE\)
be a locally Lipshitzian vector field
on an open set \(\OpenE\subset\RR^\Dim\)
with associated (local) flow %\(\VFlow{\VFieldE}{\FlowVar}{\pEuc}\)
\begin{equation*}
	\RR\times\OpenE\supset\DFlowX\VFieldE\ni(\FlowVar,\pEuc)\mapsto\VFlow{\VFieldE}{\FlowVar}{\pEuc}\in\OpenE
\end{equation*}
and let
\(\funcA,\funcB\in\LOneLoc{\OpenE}\) be given
and suppose that \(\funcB=L_\VFieldE\funcA\).

Let \(\SubOpenE\subset\OpenE\) an open set and \(\FlowVar_0>0\) be given and
assume that \([0,\FlowVar_0]\times\SubOpenE\subset\DFlowX\VField\).
Then for each \(\FlowVar\in]0,\FlowVar_0[\) the equality
%and \(0<\FlowVar\leq\FlowVar_0\)
%then for almost all \(\pEuc\in\SubOpenE\)
%(which may depend on \(\FlowVar	\))
\begin{equation*}
	\DeltaR\VField\FlowVar\funcA(\pEuc)%\funcA\bigl(\VFlow{\VField}{\FlowVar}{\pEuc}\bigr)-\funcA(\pEuc)
	=\MeanOp{\VField}{\FlowVar}{\funcB}(\pEuc)
	=\dfrac{1}{\FlowVar}\int_0^\FlowVar\funcB\bigl(\VFlow{\VField}{\FlowVarB}{\pEuc}\bigr)\, d\FlowVarB
\end{equation*}
%where
%\begin{equation*}
%	\DeltaR\VField\FlowVar\funcA(\pEuc)
%	\definedby
%	\dfrac{\funcA\bigl(\VFlow{\VField}{\FlowVar}{\pEuc}\bigr)-\funcA(\pEuc)}{\FlowVar}
%\end{equation*}
holds for almost all \(\pEuc\in\SubOpenE\).
\end{proposition}

{
\def\hVar{h}
\def\hIVar{u}
\def\SmallA#1{I_{#1}}
\def\SmallB#1{J_{#1}}
\begin{proof}
By a standard localization argument we may assume that \(\overline{\SubOpenE}\) is a compact subset of \(\OpenE\)
and that the vector field \(\VFieldE\) has compact support in \(\OpenE\) and hence is complete.
 
Let \(\FlowVar\in]0,\FlowVar_0[\) be fixed.
Since \(\funcB=L_\VFieldE\funcA\) it follows that for almost all \(\hIVar>0\) small enough and
almost \(\pEuc\in\SubOpenE\)
\begin{equation*}
	\funcA\bigl(\VFlow{\VField}{\FlowVar+\hIVar}{\pEuc}\bigr)-\funcA(\pEuc)
	=\int_0^{\FlowVar+\hIVar}\funcB\bigl(\VFlow{\VField}{\FlowVarB}{\pEuc}\bigr)\, d\FlowVarB
\end{equation*}
and hence, for each \(\hVar>0\) small enought,
\begin{equation*}
	\dfrac{1}{\hVar}\int_{0}^\hVar\bigl(\funcA\bigl(\VFlow{\VField}{\FlowVar+\hIVar}{\pEuc}\bigr)-\funcA(\pEuc)\bigr)\, d\hIVar
	=\dfrac{1}{\hVar}\int_{\hIVar=0}^\hVar\int_{\FlowVarB=0}^{\FlowVar+\hIVar}\funcB\bigl(\VFlow{\VField}{\FlowVarB}{\pEuc}\bigr)\, d\FlowVarB\, d\hIVar.
\end{equation*}
After some straightforward computations we obtain
\begin{equation*}
	\DeltaR\VField\FlowVar\funcA(\pEuc)=
	\MeanOp{\VField}{\FlowVar}{\funcB}(\pEuc)
	%\dfrac{1}{\FlowVar}\int_0^\FlowVar\funcB\bigl(\VFlow{\VField}{\FlowVarB}{\pEuc}\bigr)\, d\FlowVarB\
	-\SmallA{\hVar}(\pEuc)
	+\SmallB{\hVar}(\pEuc)
\end{equation*}
where
\begin{equation*}
	\SmallA{\hVar}(\pEuc)=\dfrac{1}{\FlowVar\,\hVar}\int_{\hIVar=0}^\hVar
		\Bigl(\funcA\circ\VFlowB{\VField}{\FlowVar}{\VFlow{\VField}{\hIVar}{\pEuc}}-\funcA\circ\VFlow{\VField}{\FlowVar}{\pEuc}\Bigr)\, d\hIVar
\end{equation*}
and
\begin{equation*}
	\SmallB{\hVar}(\pEuc)=
	\dfrac{1}{\FlowVar\,\hVar}\int_{\hIVar=0}^\hVar\int_{\FlowVarB=0}^{\hIVar}\funcB\circ\VFlowB{\VField}{\FlowVar}{\VFlow{\VField}{\FlowVarB}{\pEuc}}\, d\FlowVarB\, d\hIVar.
\end{equation*}
We complete the proof observing that the functions \(\funcA\circ\VFlowW{\VField}{\FlowVar}\) and \(\funcB\circ\VFlowW{\VField}{\FlowVar}\)
belong to \(\LOneLoc{\OpenE}\)
and hence Proposition \ref{stm::LipVFLebesgueP} easily implies that
\begin{equation*}
	\lim_{\hVar\to0^+}\SmallA{\hVar}(\pEuc)=\lim_{\hVar\to0^+}\SmallB{\hVar}(\pEuc)=0
\end{equation*}
for almost all \(\pEuc\in\SubOpenE\).

\end{proof}
}

A similar argument yields:

\begin{proposition}\label{stm::AlmostUGToUG}
Let \(\VFieldE\)
be a locally Lipshitzian vector field
on an open set \(\OpenE\subset\RR^\Dim\)
with associated (local) flow %\(\VFlow{\VFieldE}{\FlowVar}{\pEuc}\)
\begin{equation*}
	\RR\times\OpenE\supset\DFlowX\VFieldE\ni(\FlowVar,\pEuc)\mapsto\VFlow{\VFieldE}{\FlowVar}{\pEuc}\in\OpenE
\end{equation*}
and let
\(\funcA,\UGradA\in\LOneLoc{\OpenE}\) be given
and suppose that \(\UGradA\) is an upper gradient for the function \(\funcA\in\LOneLoc{\OpenE}\)
along the vector field \(\VFieldE\).

Let \(\SubOpenE\subset\OpenE\) an open set and \(\FlowVar_0>0\) be given and
assume that \([0,\FlowVar_0]\times\SubOpenE\subset\DFlowX\VField\).
Then for each \(\FlowVar\in]0,\FlowVar_0[\) the inequality
%and \(0<\FlowVar\leq\FlowVar_0\)
%then for almost all \(\pEuc\in\SubOpenE\)
%(which may depend on \(\FlowVar	\))
\begin{equation*}
	\DeltaR\VField\FlowVar\funcA(\pEuc)%\funcA\bigl(\VFlow{\VField}{\FlowVar}{\pEuc}\bigr)-\funcA(\pEuc)
	\leq\MeanOp{\VField}{\FlowVar}{\UGradA}(\pEuc)
	=\dfrac{1}{\FlowVar}\int_0^\FlowVar\UGradA\bigl(\VFlow{\VField}{\FlowVarB}{\pEuc}\bigr)\, d\FlowVarB
\end{equation*}
%where
%\begin{equation*}
%	\DeltaR\VField\FlowVar\funcA(\pEuc)
%	\definedby
%	\dfrac{\funcA\bigl(\VFlow{\VField}{\FlowVar}{\pEuc}\bigr)-\funcA(\pEuc)}{\FlowVar}
%\end{equation*}
holds for almost all \(\pEuc\in\SubOpenE\).
\end{proposition}

}%end unit

%\section{\label{section:VFOperators}Basic Operators}
\subsection{Difference quotients}
% !TeX encoding = UTF-8
% !TeX spellcheck = en_GB

{%begin unit
\def\fOne{\funcA}
\def\fTwo{\funcTest}

In this section we prove that
for each \(\funcA\in\LOneLoc{\OpenE}\)
%further consequences of Proposition \ref{stm::VFJacToDiv}.
the difference quotients
\(%\begin{equation*}
	\DeltaR\VField\FlowVar\funcA(\pEuc)
%	\definedby
%	\dfrac{\funcA\bigl(\VFlow{\VField}{\FlowVar}{\pEuc}\bigr)-\funcA(\pEuc)}{\FlowVar}
\) %\end{equation*}
always converge as a distribution when \(\FlowVar\to0\):
see Proposition \ref{stm::DQLim}.
%associated to a (locally) integrable function \(\funcA\) on an \(\OpenE\).

\begin{lemma}
Let \(\VField\)
be a complete Lipshitzian vector field
on an open set \(\OpenE\subset\RR^\Dim\).
%with associated flow
%\begin{equation*}
%	\RR\times\OpenE\ni(\FlowVar,\pEuc)\mapsto\VFlow{\VField}{\FlowVar}{\pEuc}\in\OpenE.
%\end{equation*}
If \(\fOne\in L^1(\OpenE)\) then
\begin{equation*}
	\lim_{\FlowVar\to0}
		\int_\OpenE\DeltaR\VField\FlowVar\funcA(\pEuc)
			%\dfrac{\fOne\bigl(\VFlow{\VFieldE}{\FlowVar}{\pEuc}\bigr)
			%-\fOne(\pEuc)}{\FlowVar}\,d\pEuc
	=-\int_\OpenE\fOne(\pEuc)\Divergence\VFieldE(\pEuc)\,d\pEuc.
\end{equation*}
\end{lemma}

\begin{proof}
Let \(\FlowVar\in\RR\), \(\FlowVar\neq0\),	 be given.
The change of variables
\begin{equation*}
	\pEuc\mapsto\VFlow{\VFieldE}{-\FlowVar}{\pEuc}
\end{equation*}
yields
\begin{equation*}
	\int_\OpenE\fOne\bigl(\VFlow{\VFieldE}{\FlowVar}{\pEuc}\bigr)\,d\pEuc
	=\int_\OpenE\fOne(\pEuc)\VFJac\VFieldE{-\FlowVar}\pEuc\,d\pEuc
\end{equation*}
and hence
\begin{equation*}
\begin{split}
		\int_\OpenE\DeltaR\VField\FlowVar\funcA(\pEuc)
%			\dfrac{\fOne\bigl(\VFlow{\VFieldE}{\FlowVar}{\pEuc}\bigr)
%				-\fOne(\pEuc)}{\FlowVar}\,d\pEuc
		&=
		\dfrac{1}{\FlowVar}
			\int_\OpenE\fOne\bigl(\VFlow{\VFieldE}{\FlowVar}{\pEuc}\bigr)\,d\pEuc
		-\dfrac{1}{\FlowVar}
			\int_\OpenE\fOne(\pEuc)\,d\pEuc
		\\
		&=
		\dfrac{1}{\FlowVar}
			\int_\OpenE\fOne(\pEuc)\VFJac\VFieldE{-\FlowVar}\pEuc\,d\pEuc
		-\dfrac{1}{\FlowVar}
			\int_\OpenE\fOne(\pEuc)\,d\pEuc
		\\
		&=
		-\int_\OpenE\fOne(\pEuc)
			\dfrac{\VFJac\VFieldE{-\FlowVar}\pEuc-1}{-\FlowVar}\,d\pEuc.		
\end{split}
\end{equation*}
Proposition \ref{stm::VFJacToDiv} implies then
\begin{equation*}
\begin{split}
		\lim_{\FlowVar\to0}
		\int_\OpenE\DeltaR\VField\FlowVar\funcA(\pEuc)
%			\dfrac{\fOne\bigl(\VFlow{\VFieldE}{\FlowVar}{\pEuc}\bigr)
%			-\fOne(\pEuc)}{\FlowVar}\,d\pEuc
		&=-\lim_{\FlowVar\to0}
			\int_\OpenE\fOne(\pEuc)
			\dfrac{\VFJac\VFieldE{-\FlowVar}\pEuc-1}{-\FlowVar}\,d\pEuc
		\\&=-\int_\OpenE\fOne(\pEuc)\Divergence\VFieldE(\pEuc)\,d\pEuc,
\end{split}
\end{equation*}
as required.

\end{proof}

Let denote by \(\LipC\OpenE\) the space of
compactly supported lipshitzian real functions on \(\OpenE\)

\begin{proposition}\label{stm::DQLim}
Let \(\VField\)
be a locally lipshitzian vector field
on an open set \(\OpenE\subset\RR^\Dim\)
%with associated flow
%\begin{equation*}
%	\RR\times\OpenE\ni(\FlowVar,\pEuc)\mapsto\VFlow{\VField}{\FlowVar}{\pEuc}\in\OpenE.
%\end{equation*}
If \(\fOne\in\LOneLoc{\OpenE}\)
%and \(\fTwo\in C_0^1(\OpenE)\)
and \(\fTwo\in\LipC{\OpenE}\)
then
\begin{equation*}
	\lim_{\FlowVar\to0}
	\int_\OpenE\DeltaR\VField\FlowVar\funcA(\pEuc)\fTwo(\pEuc)\,d\pEuc
	%\dfrac{\fOne\bigl(\VFlow{\VFieldE}{\FlowVar}{\pEuc}\bigr)
		%	-\fOne(\pEuc)}{\FlowVar}&\fTwo(\pEuc)\,d\pEuc
	=-\int_\OpenE\fOne(\pEuc)\VFieldE\fTwo(\pEuc)\,d\pEuc
	-\int_\OpenE\fOne(\pEuc)\fTwo(\pEuc)\Divergence\VFieldE(\pEuc)\,d\pEuc.
\end{equation*}
%\begin{equation*}
%\begin{split}
%	\lim_{\FlowVar\to0}
%	\int_\OpenE\DeltaR\VField\FlowVar\funcA(\pEuc)
%	%\dfrac{\fOne\bigl(\VFlow{\VFieldE}{\FlowVar}{\pEuc}\bigr)
%	%	-\fOne(\pEuc)}{\FlowVar}&\fTwo(\pEuc)\,d\pEuc
%	\\=-\int_\OpenE\fOne(\pEuc)&\VFieldE\fTwo(\pEuc)\,d\pEuc
%	-\int_\OpenE\fOne(\pEuc)\fTwo(\pEuc)\Divergence\VFieldE(\pEuc)\,d\pEuc.
%\end{split}
%\end{equation*}
\end{proposition}

\begin{proof}
We may replace \(\VFieldE\)
with a vector field
which coincides with \(\VFieldE\)
in a neighbourhood of the support of the function \(\fTwo\)
and
having compact support in \(\OpenE\).
The same can be done for the function \(\fOne\)
%(see e.g., Proposition \ref{stm::LocalizationPrinciple})
.

Doing so we can assume that
the vector field \(\VFieldE\) is complete
and that \(\fOne\in L^1(\OpenE)\).

Let \(\FlowVar\in\RR\), \(\FlowVar\neq0\), be given.
The same arguments used in the proof of the previous lemma
and some straightforward algebraic manipulations
yield the identity
\begin{equation*}
\begin{split}
	\int_\OpenE\DeltaR\VField\FlowVar\funcA(\pEuc)\fTwo(\pEuc)\,&d\pEuc
	%\dfrac{\fOne\bigl(\VFlow{\VFieldE}{\FlowVar}{\pEuc}\bigr)
	%-\fOne(\pEuc)}{\FlowVar}&\fTwo(\pEuc)\,d\pEuc
	\\
	=
	-\int_\OpenE\fOne(\pEuc)&\DeltaR\VField{-\FlowVar}\fTwo(\pEuc)
		%\dfrac{\fTwo\bigl(\VFlow{\VFieldE}{-\FlowVar}{\pEuc}\bigr)
		%	-\fTwo(\pEuc)}{-\FlowVar}
	\VFJac\VFieldE{-\FlowVar}\pEuc\,d\pEuc
	-\int_\OpenE\fOne(\pEuc)\fTwo(\pEuc)
	\dfrac{\VFJac\VFieldE{-\FlowVar}\pEuc-1}{-\FlowVar}\,d\pEuc.		
\end{split}
\end{equation*}

We have
\begin{equation*}
	\LimVFZ\VFJac\VFieldE{-\FlowVar}\pEuc
	=\LimVFZ\VFJac\VFieldE{\FlowVar}\pEuc
	=1
\end{equation*}
being
\(%\begin{equation*}
%\abs{\VFJac\VFieldE{-\FlowVar}\pEuc-1}
\abs{\VFJac\VFieldE{\FlowVar}\pEuc-1}
= O(\abs{\FlowVar})
\) %\end{equation*}
by Proposition \ref{stm::VFJacToDiv}.
The classical Rademacher theorem implies that
for almost all \(\pEuc\in\OpenE\)
\begin{equation*}
	\LimVFZ\DeltaR\VField{-\FlowVar}\fTwo(\pEuc)%\dfrac{\fTwo\bigl(\VFlow{\VFieldE}{-\FlowVar}{\pEuc}\bigr)
		%-\fTwo(\pEuc)}{-\FlowVar}
	=\LimVFZ\DeltaR\VField\FlowVar\fTwo(\pEuc)%\dfrac{\fTwo\bigl(\VFlow{\VFieldE}{\FlowVar}{\pEuc}\bigr)
		%-\fTwo(\pEuc)}{\FlowVar}
	=\VFieldE\fTwo(\pEuc)
\end{equation*}
since \(\fTwo\in\LipC{\OpenE}\).%\(\fTwo\in C_0^1(\OpenE)\).

Observe now that for any fixed \(\FlowVar_0>0\) the functions
\begin{equation*}
	\OpenE\ni\pEuc\mapsto\DeltaR\VField\FlowVar\fTwo(\pEuc)
		%\dfrac{\fTwo\bigl(\VFlow{\VFieldE}{\FlowVar}{\pEuc}\bigr)
		%-\fTwo(\pEuc)}{\FlowVar}
	\VFJac\VFieldE{-\FlowVar}\pEuc
\end{equation*}
are uniformly bounded when \(-\FlowVar_0\leq\FlowVar\leq\FlowVar_0\)
%Both limits holds uniformly with respect to \(\pEuc\in\OpenE\) 
and hence the Lebesgue theorem on the dominate convergence yields
\begin{equation*}
	\LimVFZ\int_\OpenE\fOne(\pEuc)\DeltaR\VField{-\FlowVar}\fTwo(\pEuc)
	%\dfrac{\fTwo\bigl(\VFlow{\VFieldE}{-\FlowVar}{\pEuc}\bigr)
	%	-\fTwo(\pEuc)}{-\FlowVar}
	\VFJac\VFieldE{-\FlowVar}\pEuc\,d\pEuc
	=\int_\OpenE\fOne(\pEuc)\VFieldE\fTwo(\pEuc)\,d\pEuc
\end{equation*}
Proposition \ref{stm::VFJacToDiv}
(and the Lebesgue theorem on the dominate convergence)
again yields
\begin{equation*}
	\LimVFZ\int_\OpenE\fOne(\pEuc)\fTwo(\pEuc)
	\dfrac{\VFJac\VFieldE{-\FlowVar}\pEuc-1}{-\FlowVar}\,d\pEuc
	=\int_\OpenE\fOne(\pEuc)\fTwo(\pEuc)\Divergence\VFieldE(\pEuc)\,d\pEuc
\end{equation*}
and hence
\begin{equation*}
\begin{split}
	\LimVFZ&\int_\OpenE\DeltaR\VField\FlowVar\funcA(\pEuc)\fTwo(\pEuc)\,d\pEuc
	%=\LimVFZ\int_\OpenE&
	%\dfrac{\fOne\bigl(\VFlow{\VFieldE}{\FlowVar}{\pEuc}\bigr)
	%	-\fOne(\pEuc)}{\FlowVar}\fTwo(\pEuc)\,d\pEuc
	\\
	&=
	-\LimVFZ\int_\OpenE\fOne(\pEuc)\DeltaR\VField{-\FlowVar}\fTwo(\pEuc)
	%\dfrac{\fTwo\bigl(\VFlow{\VFieldE}{-\FlowVar}{\pEuc}\bigr)
	%	-\fTwo(\pEuc)}{-\FlowVar}
	\VFJac\VFieldE{-\FlowVar}\pEuc\,d\pEuc
	%\\
	%&
	-\LimVFZ\int_\OpenE\fOne(\pEuc)\fTwo(\pEuc)
	\dfrac{\VFJac\VFieldE{-\FlowVar}\pEuc-1}{-\FlowVar}\,d\pEuc
	\\
	&=-\int_\OpenE\fOne(\pEuc)\VFieldE\fTwo(\pEuc)\,d\pEuc
	-\int_\OpenE\fOne(\pEuc)\fTwo(\pEuc)\Divergence\VFieldE(\pEuc)\,d\pEuc,
\end{split}
\end{equation*}
%\begin{equation*}
%\begin{split}
%	\int_\OpenE\DeltaR\VField\FlowVar\funcA(\pEuc)\fTwo(\pEuc)\,d\pEuc=
%	\LimVFZ&\int_\OpenE
%	\dfrac{\fOne\bigl(\VFlow{\VFieldE}{\FlowVar}{\pEuc}\bigr)
%		-\fOne(\pEuc)}{\FlowVar}\fTwo(\pEuc)\,d\pEuc
%	\\
%	=
%	-\LimVFZ\int_\OpenE\fOne(\pEuc)&
%	\dfrac{\fTwo\bigl(\VFlow{\VFieldE}{-\FlowVar}{\pEuc}\bigr)
%		-\fTwo(\pEuc)}{-\FlowVar}\VFJac\VFieldE{-\FlowVar}\pEuc\,d\pEuc
%	\\
%	-\LimVFZ\int_\OpenE&\fOne(\pEuc)\fTwo(\pEuc)
%	\dfrac{\VFJac\VFieldE{-\FlowVar}\pEuc-1}{-\FlowVar}\,d\pEuc
%	\\
%	=-\int_\OpenE\fOne(\pEuc)&\VFieldE\fTwo(\pEuc)\,d\pEuc
%	-\int_\OpenE\fOne(\pEuc)\fTwo(\pEuc)\Divergence\VFieldE(\pEuc)\,d\pEuc,
%\end{split}
%\end{equation*}
as desired.

\end{proof}

}%end unit

\section{\label{section:VFComplete}Complete Lipshitzian vector fields}
% !TeX encoding = UTF-8
% !TeX spellcheck = en_GB
{%begin unit
In this section,
we assume that
\(\VField\) is a complete Lipshitzian vector field
on an open set \(\OpenE\subset\RR^\Dim\)
with Lipshitz constant \(\LipVF\)
and
associated flow
\begin{equation*}
	\RR\times\OpenE\ni(\FlowVar,\pEuc)\mapsto\VFlow{\VField}{\FlowVar}{\pEuc}\in\OpenE.
\end{equation*}

%\(\GPFlow{\VFieldE}{\FlowVar}\) is the induced
%\(c_0\) one parameter group on the Banach space \(L^1(\OpenE)\),
%\begin{equation*}
%	L^1(\OpenE)\ni\funcA\mapsto\GPFlow\VField\FlowVar\funcA=\funcA\circ\VFlowW{\VField}{\FlowVar}\in L^1(\OpenE)
%\end{equation*}
%and \(\GFZ:\Dom(\GFZ)\to L^1(\OpenE)\) is its infinitesimal
%generator.

}%end unit

\subsection{The mean operator}\label{ssec::MeanOp}
% !TeX encoding = UTF-8
% !TeX spellcheck = en_GB

{%begin unit
%Now we study the action of the operator
%\(\MeanOp{\VFieldE}{\FlowVar}\) on \(L^1(\OpenE)\).
%In this section
%\(\OpenE\sset\RR^\Dim\) is an open set
%and
%\(\VFieldE\) is a complete lipshitzian vector field
%with Lipshitz constant \(\LipVF\)
%and
%associated flow \(\VFlow{\VFieldE}{\FlowVar}{\cdot}\).
%
%We define the \emph{mean operator} associated to the vector field \(\VFieldE\)
%as the operator defined for each \(\funcA\in\LOneLoc{\OpenE}\) by the formula
%\begin{equation}\label{stm::MeanOpDef}
%	\MeanOp{\VFieldE}{\FlowVar}{\funcA}(\pEuc)=\dfrac{1}{\FlowVar}
%	\int_0^\FlowVar\funcA\bigl(\VFlow{\VFieldE}{\FlowVarB}{\pEuc}\bigr)\,d\FlowVarB
%\end{equation}
%when \(\FlowVar\neq0\);
%we also set \(\MeanOp{\VFieldE}{0}{\funcA}=\funcA\).
%Recalling that we assume that
%Sthe vector field is lipshitzian and complete then we have:
\begin{lemma}\label{stm::GFDefined}
If \(\funcA\in L^1(\OpenE)\)
and \(\FlowVar\in\RR\)
then 
\(\funcA\circ\VFlowW{\VField}{\FlowVar}\in L^1(\OpenE)\)
and
\begin{equation*}\label{stm::VFLebesgueP}
	\norm{\funcA\circ\VFlowW{\VField}{\FlowVar}}_{L^1(\OpenE)}
	\leq
	e^{\Dim\LipVF\abs{\FlowVar}}\norm{\funcA}_{L^1(\OpenE)}.
\end{equation*}
\end{lemma}

\begin{proof}
Using the change of variable
\(\pEuc\mapsto\VFlow{\VField}{-\FlowVar}{\pEuc}\)
and the estimate \eqref{stm::VFJExpEstimate}
we obtain
\begin{equation*}
	\norm{\funcA\circ\VFlowW{\VField}{\FlowVar}}_{L^1(\OpenE)}
	=
	\int_\OpenE\abs{\funcA\bigl(\VFlow{\VField}{\FlowVar}{\pEuc}\bigr)}\,d\pEuc
	=
	\int_\OpenE\abs{\funcA(\pEuc)}\VFJac\VFieldE{-\FlowVar}\pEuc\,d\pEuc
	\leq e^{\Dim\LipVF\abs{\FlowVar}}\norm{\funcA}_{L^1(\OpenE)}.
\end{equation*}

\end{proof}

\begin{lemma}\label{stm::FlowIsCZero}
For each \(\funcA\in L^1(\OpenE)\) we have
\begin{equation}\label{eq::FlowIsCZero}
	\lim_{\FlowVar\to0}\norm{\funcA\circ\VFlowW{\VField}{\FlowVar}-\funcA}_{L^1(\OpenE)}
	=0.
\end{equation}
\end{lemma}

\begin{proof}
It is easy to show that \eqref{eq::FlowIsCZero} hold if \(\funcA\in C_0(\OpenE)\).
An elementary argument using
the density of \(C_0(\OpenE)\) in \(L^1(\OpenE)\) and
the estimate given in Lemma \ref{stm::GFDefined}
imply the desired assertion.

\end{proof}

\begin{proposition}\label{stm::MeanOpMain}
If \(\funcA\in L^1(\OpenE)\) then
%the formula \eqref{stm::MeanOpDef} defines a continuous linear operator
\(%\begin{equation*}
	\MeanOp{\VFieldE}{\FlowVar}\funcA\in L^1(\OpenE)%\to L^1(\OpenE)
\) % \end{equation*}
for each \(\FlowVar\in\RR\)
and 
%for each \(\funcA\in L^1(\OpenE)\)
the map
\begin{equation*}
	\RR\ni\FlowVar\mapsto\MeanOp{\VFieldE}{\FlowVar}\funcA\in L^1(\OpenE)
\end{equation*}
is continuous when \(L^1(\OpenE)\)
is endowed with the topology induced by its standard norm.
\end{proposition}

\begin{proof}{
\def\fva{{\FlowVar_1}}
\def\fvb{{\FlowVar_2}}
\def\LocalOP#1{K_{#1}}
The operator \(\MeanOp{\VField}{0}\)
is the identity map on \(L^1(\OpenE)\) and hence is obviously bounded.

Let \(\FlowVar>0\) and \(\funcA\in L^1(\OpenE)\).
%Proposition \ref{stm::GPversusSGP} allow us to suppose \(\FlowVar\geq0\);
%the change of variables
%\begin{equation*}
%	\pEuc\mapsto\VFlow{\VField}{-\FlowVar}{\pEuc}
%\end{equation*}
Lemma \ref{stm::GFDefined}
%and (\rmnum{\stmJBouned}) of Proposition \ref{stm::VFJacToDiv}
yields
\begin{equation*}
\begin{split}
	\norm{\MeanOp{\VField}{\FlowVar}{\funcA}}_{L^1(\OpenE)}
		&=\dfrac{1}{\FlowVar}\int_\OpenE
			\abs{\int_0^\FlowVar\funcA\bigl(\VFlow{\VField}{\FlowVarB}{\pEuc}\bigr)\,d\FlowVarB}d\pEuc
		\\
		 &\leq\dfrac{1}{\FlowVar}\int_\OpenE\int_0^\FlowVar
			\abs{\funcA\bigl(\VFlow{\VField}{\FlowVarB}{\pEuc}\bigr)}d\FlowVarB\,d\pEuc
		 =\dfrac{1}{\FlowVar}\int_0^\FlowVar\int_\OpenE
			\abs{\funcA\bigl(\VFlow{\VField}{\FlowVarB}{\pEuc}\bigr)}d\pEuc\,d\FlowVarB
		\\
		% =\dfrac{1}{\FlowVar}\int_0^\FlowVar\int_\OpenE
		%	\abs{\funcA(\pEuc)}\VFJac\VFieldE{-\FlowVar}{\pEuc}\,d\pEuc\,d\FlowVarB
		&=\dfrac{1}{\FlowVar}\int_0^\FlowVar\norm{\funcA\circ\VFlowMap{\VFieldE}{\FlowVarB}}_{L^1(\OpenE)}d\FlowVarB
		 \leq\dfrac{1}{\FlowVar}\int_0^\FlowVar e^{\Dim\LipVF\FlowVarB}\norm{\funcA}_{L^1(\OpenE)}d\FlowVarB
		\\
		&=\dfrac{e^{\Dim\LipVF\FlowVar}-1}{\Dim\LipVF\FlowVar}\norm{\funcA}_{L^1(\OpenE)}
%		&\leq e^{\Dim\LipVF\abs{\FlowVar}}\dfrac{1}{\FlowVar}\int_0^\FlowVar\int_\OpenE
%			\abs{\funcA(\pEuc)}d\pEuc\,d\FlowVarB
%		 =e^{\Dim\LipVF\abs{\FlowVar}}\norm{\funcA}_{L^(\OpenE)}.
\end{split}
\end{equation*}
Since \(\funcA\in L^1(\OpenE)\) is arbitrary this proves the
boundedness of the operator \(\MeanOp{\VFieldE}{\FlowVar}\) on \(L^1(\OpenE)\).

For each \(\FlowVar>0\) define
\begin{equation*}
	\LocalOP\FlowVar\funcA(\pEuc)=
		\int_0^\FlowVar\funcA\bigl(\VFlow{\VFieldE}{\FlowVarB}{\pEuc}\bigr)\,d\FlowVarB.
\end{equation*}
When \(0<\fva<\fvb\) we have
\begin{equation*}
\begin{split}
	\LocalOP{\fvb}\funcA(\pEuc)-\LocalOP{\fva}\funcA(\pEuc)
	=\int_\fva^\fvb\funcA\bigl(\VFlow{\VFieldE}{\FlowVarB}{\pEuc}\bigr)\,d\FlowVarB
	%-\int_0^\fva\funcA\bigl(\VFlow{\VFieldE}{\FlowVarB}{\pEuc}\bigr)\,d\FlowVarB
\end{split}
\end{equation*}
and hence, by Lemma \ref{stm::GFDefined} again,
\begin{equation*}
\begin{split}
	\norm{\LocalOP{\fvb}\funcA-\LocalOP{\fva}\funcA}_{L^1(\OpenE)}
	&=\int_\OpenE\abs{\int_\fva^\fvb\funcA\bigl(\VFlow{\VFieldE}{\FlowVarB}{\pEuc}\bigr)\,d\FlowVarB}d\pEuc
	\\
	&\leq\int_\OpenE\int_\fva^\fvb\abs{\funcA\bigl(\VFlow{\VFieldE}{\FlowVarB}{\pEuc}\bigr)}d\FlowVarB\,d\pEuc
	\\
	&=\int_\fva^\fvb\int_\OpenE\abs{\funcA\bigl(\VFlow{\VFieldE}{\FlowVarB}{\pEuc}\bigr)}d\pEuc\,d\FlowVarB
	\\
	&=\int_\fva^\fvb%\int_\OpenE
	%\abs{\funcA(\pEuc)}\VFJac\VFieldE{-\FlowVar}{\pEuc}\,d\pEuc\,d\FlowVarB
	\norm{\funcA\circ\VFlowMap{\VFieldE}{\FlowVarB}}_{L^1(\OpenE)}d\FlowVarB
	\\
	&\leq\int_\fva^\fvb e^{\Dim\LipVF\FlowVarB}\norm{\funcA}_{L^1(\OpenE)}d\FlowVarB
	%&\leq\int_\fva^\fvb e^{\Dim\LipVF\FlowVarB}\,d\FlowVarB\int_\OpenE\abs{\funcA(\pEuc)}d\pEuc
	\\
	&=\dfrac{e^{\Dim\LipVF\fvb}-e^{\Dim\LipVF\fva}}{\Dim\LipVF}\norm{\funcA}_{L^1(\OpenE)}.
\end{split}
\end{equation*}
Since \(\MeanOp{\VFieldE}{\FlowVar}{\funcA}=\FlowVar^{-1}\LocalOP{\FlowVar}\funcA\)
it follows that the function
\(\FlowVar\mapsto\MeanOp{\VFieldE}{\FlowVar}{\funcA}\)
is continuous when \(\FlowVar>0\).

Let us now prove the right continuity of
\(\FlowVar\mapsto\MeanOp{\VFieldE}{\FlowVar}{\funcA}\)
at \(\FlowVar=0\).

For each \(\FlowVar>0\) we have
\begin{equation*}
\begin{split}
	\norm{\MeanOp{\VFieldE}{\FlowVar}{\funcA}-\funcA}_{L^1(\OpenE)}
	&=\int_\OpenE
	\abs{\dfrac{1}{\FlowVar}\int_0^\FlowVar\funcA\bigl(\VFlow{\VField}{\FlowVarB}{\pEuc}\bigr)\,d\FlowVarB
		-\funcA(\pEuc)}d\pEuc
	\\
	&=\int_\OpenE
	\abs{\dfrac{1}{\FlowVar}\int_0^\FlowVar\Bigl(\funcA\bigl(\VFlow{\VField}{\FlowVarB}{\pEuc}\bigr)-\funcA(\pEuc)\Bigr)
		\,d\FlowVarB}d\pEuc
	\\
	&\leq\dfrac{1}{\FlowVar}\int_\OpenE
		\int_0^\FlowVar\abs{\funcA\bigl(\VFlow{\VField}{\FlowVarB}{\pEuc}\bigr)-\funcA(\pEuc)}
		d\FlowVarB\,d\pEuc
	\\
	&=\dfrac{1}{\FlowVar}\int_0^\FlowVar\int_\OpenE
	\abs{\funcA\bigl(\VFlow{\VField}{\FlowVarB}{\pEuc}\bigr)-\funcA(\pEuc)}
		d\pEuc\,d\FlowVarB
	\\
	&=\dfrac{1}{\FlowVar}\int_0^\FlowVar
		\norm{\funcA\circ\VFlowW{\VField}{\FlowVarB}-\funcA}_{L^1(\OpenE)}d\FlowVarB
\end{split}
\end{equation*}
Lemma \ref{stm::FlowIsCZero} says that
\begin{equation*}
	\lim_{\FlowVar\to0}\norm{\funcA\circ\VFlowW{\VField}{\FlowVar}-\funcA}_{L^1(\OpenE)}
	=0.
\end{equation*}
and hence we obtain 
\begin{equation*}
	\lim_{\FlowVar\to0^+}\norm{\MeanOp{\VFieldE}{\FlowVar}{\funcA}-\funcA}_{L^1(\OpenE)}=0,
\end{equation*}
as desired.

The same argument applied to the vector field \(-\VFieldE\)
completes the proof.

}\end{proof}

We end this section proving a curious relation
between the mean operator and the different quotient
which will be used in the proof of Theorem \ref{stm::LOne::ZEqLEqX}.

\begin{proposition}
Let \(\FlowVar, \FlowVarM\in\RR\setminus\{0\}\) be given.
Then the difference quotient \(\DeltaR\VField\FlowVar{}\)
and the mean operator \(\MeanOp\VField\FlowVar{}\)
satisfy the relation
\begin{equation*}
	\MeanOp\VField\FlowVarM{}\DeltaR\VField\FlowVar{}
	=
	\MeanOp\VField\FlowVar{}\DeltaR\VField\FlowVarM{}
\end{equation*} 
\end{proposition} 

\begin{proof}
Let \(\funcA\in L^1(\OpenE)\) be given.
The group identities like
\begin{equation*}
	\VFlow{\VFieldE}{\FlowVar}{\VFlow{\VFieldE}{\FlowVarM}{\pEuc}}
	=\VFlow{\VFieldE}{\FlowVar+\FlowVarM}{\pEuc}
\end{equation*}
and some straightforward computations yield
\begin{equation*}
	\MeanOp\VField\FlowVarM{}\DeltaR\VField\FlowVar\funcA(\pEuc)
	=\dfrac{1}{\FlowVarM\FlowVar}
	\left(\int_\FlowVar^{\FlowVar+\FlowVarM}\funcA\bigl(\VFlow{\VField}{\FlowVarB}{\pEuc}\bigr)d\FlowVarB
		-\int_0^{\FlowVarM}\funcA\bigl(\VFlow{\VField}{\FlowVarB}{\pEuc}\bigr)d\FlowVarB
	\right)
\end{equation*}
and
\begin{equation*}
	\MeanOp\VField\FlowVar{}\DeltaR\VField\FlowVarM\funcA(\pEuc)
	=\dfrac{1}{\FlowVarM\FlowVar}
	\left(\int_\FlowVarM^{\FlowVar+\FlowVarM}\funcA\bigl(\VFlow{\VField}{\FlowVarB}{\pEuc}\bigr)d\FlowVarB
	-\int_0^{\FlowVar}\funcA\bigl(\VFlow{\VField}{\FlowVarB}{\pEuc}\bigr)d\FlowVarB
	\right).
\end{equation*}
We conclude observing that the (Lebesgue) integral on the real line satisfies
\begin{equation*}
	\int_\FlowVar^{\FlowVar+\FlowVarM}-\int_0^{\FlowVarM}
	=\int_\FlowVarM^{\FlowVar+\FlowVarM}-\int_0^{\FlowVar}.
\end{equation*}

 \end{proof}

}%end unit

\subsection{The associated one parameter group}
% !TeX encoding = UTF-8
% !TeX spellcheck = en_GB
{%begin unit
%In this section
%\(\OpenE\sset\RR^\Dim\) is an open set
%and
%\(\VFieldE\) is a complete lipshitzian vector field
%with Lipshitz constant \(\LipVF\)
%and
%associated flow \(\VFlow{\VFieldE}{\FlowVar}{\cdot}\).

\begin{theorem}\label{stm::FlowToGroup}
Let \(\VField\)  be a complete Lipshitzian vector field
on an open set \(\OpenE\subset\RR^\Dim\).
%and let %\(\VFlow{\VField}{\FlowVar}{\pEuc}\) 
%\begin{equation*}
%	\RR\times\OpenE\ni(\FlowVar,\pEuc)\mapsto\VFlow{\VField}{\FlowVar}{\pEuc}\in\OpenE
%\end{equation*}
%be the associated flow.
%Let \(1\leq p<+\infty\) be given.
Then the map %\(\VFlow{\VField}{\FlowVar}{\pEuc}\)
\begin{equation*}
	\RR\times L^1(\OpenE)\ni(\FlowVar,\funcA)\mapsto
	\GPFlow\VField\FlowVar\funcA\definedby\funcA\circ\VFlowW{\VField}{\FlowVar}\in L^1(\OpenE)
\end{equation*}
defines a jointly continuos (i. e. \(c_0\)) one parameter group
on the Banach space \(L^1(\OpenE)\)
which satisfies
\begin{equation}\label{eq::GFEstimate}
	\norm{\GPFlow\VField\FlowVar}\leq e^{\Dim\LipVF\abs{\FlowVar}}.
\end{equation}
for each \(\FlowVar\in\RR\).
\end{theorem}

\begin{proof}
Lemma \ref{stm::GFDefined} implies that \(\GPFlow\VField\FlowVar\)
is a linear operator from \(L^1(\OpenE)\) to \(L^1(\OpenE)\)
whose operator norm satisfies \eqref{eq::GFEstimate}
and
the condition \Case{1} and \Case{2} of subsection \ref{ssec::OPS}
are obviously satisfied.

Lemma \ref{stm::FlowIsCZero} implies that
\(%\begin{equation*}
\lim_{\FlowVar\to0}\norm{\GPFlow\VField\FlowVar\funcA-\funcA}_{L^1(\OpenE)}=0
\) %\end{equation*}
for each \(\funcA\in L^1(\OpenE)\)
and hence
Proposition \ref{stm::BSCtoJC} yields the desired conclusion.

\end{proof}

}%end unit

\subsection{The infinitesimal generator}
%\ifdef{\cbVOne}{
%	\input{src/VFLieDer}
%}{}
%\ifdef{\cbVTwo}{
%	\input{src/VFLieDerVTwo}
%}{}
%\ifdef{\cbVThree}{
% !TeX encoding = UTF-8
% !TeX spellcheck = en_GB
{%begin unit

Here, one of the main results of this paper.

\begin{theorem}\label{stm::LOne::ZEqLEqX}
Let \(\VField\)  be a complete Lipshitzian vector field
on an open set \(\OpenE\subset\RR^\Dim\).
Let \(\GPFlow{\VFieldE}{\FlowVar}\) be the induced
\(c_0\) one parameter group on the Banach space \(L^1(\OpenE)\)
and let \(\GFZ:\Dom(\GFZ)\to L^1(\OpenE)\) be its infinitesimal
generator.

Let \(\funcA, \funcB\in L^1(\OpenE)\) be given.
Then the following conditions are equivalent:

\begin{enumerate}[label=\cbRomanLabel]
\item\label{\stmGFX}\(\funcB=\VFieldE\funcA\), that is \(\funcB\) is
the distributional derivative of \(\funcA\);  
\item\label{\stmGFL}\(\funcB=L_\VFieldE\funcA\), that is \(\funcB\) is
the Lie derivative of \(\funcA\).
\item\label{\stmGFZ} \(\funcA\in\Dom(\GFZ)\) and \(\GFZ\funcA=\funcB\); 
\end{enumerate}
\end{theorem}

\begin{proof}{
\def\pFVar{\FlowVarM}
\def\proofTestSpace{\LipC{\OpenE}}
\def\proofTestSpace{C^\infty_0(\OpenE)}

%In order to save space we will omit in the
%formulas appearing in this proof the integration variables
To save space, in the formulas below we omit the integration variables %  in this demonstration
and the various `\(dx\)' symbols under the integral sign.

The symbol `\(\lim E\)',
where \(E\) in an expression denoting an element of \(L^1(\OpenE)\),
will always denote the limit with respect to the
topology of \(L^1(\OpenE)\) induced by its standard norm.

%We will prove that \Case{\stmGFX} and \Case{\stmGFL} are both equivalent to \Case{\stmGFZ}.

\Case{\stmGFX}\(\Longrightarrow\)\Case{\stmGFZ}:
if \(\funcB=\VFieldE\funcA\) then by Proposition \ref{stm::DQLim} %again
for each function \(\funcTest\in\proofTestSpace\) we have
\begin{equation*}
	\int_\OpenE\funcB\funcTest=
	-\int_\OpenE\funcA\VFieldE\funcTest
	-\int_\OpenE\funcA\funcTest\Divergence\VFieldE
	=\lim_{\FlowVar\to0}\int_\OpenE\DeltaR\VField\FlowVar\funcA\funcTest.
\end{equation*}
Since \(\proofTestSpace\) is dense in \(L^\infty(\OpenE)\)
with respect to the \weakstar topology then
Proposition \ref{stm::GFWStar}
%and the density of \(\proofTestSpace\) in \(L^\infty(\OpenE)\)
%with respect to the \(*-\)weak topology
implies that \(\funcA\in\Dom(\GFZ)\) and
\(\GFZ\funcA=\funcB\).

%\Case{\stmGFZ}\(\Longrightarrow\)\Case{\stmGFX}:
%since \(\funcB=\GFZ\funcA\) then, by definition of infinitesimal generator, we have
%\begin{equation*}
%	\lim_{\FlowVar\to0}\DeltaR\VField\FlowVar\funcA=g.
%\end{equation*}
%
%Then Proposition \ref{stm::DQLim} implies that
%for each \(\funcTest\in\proofTestSpace\) we have
%\begin{equation*}
%	\int_\OpenE\funcB\funcTest=
%	\lim_{\FlowVar\to0}\int_\OpenE\DeltaR\VField\FlowVar\funcA
%	\funcTest=
%	-\int_\OpenE\funcA\VFieldE\funcTest
%	-\int_\OpenE\funcA\funcTest\Divergence\VFieldE,
%\end{equation*}
%that is, by definition, \(\funcB\) is the distributional derivative of \(\funcA\).

\Case{\stmGFZ}\(\Longrightarrow\)\Case{\stmGFL}:
condition \Case{\stmGFZ} is equivalent to
\begin{equation*}
	\funcB=\lim_{\FlowVar\to0}\DeltaR{\VFieldE}{\FlowVar}{\funcA}
\end{equation*}
with respect to the topology in \(L^1(\OpenE)\) induced by its standard norm.

Using the various properties of the operator \(\MeanOp{\VFieldE}{\FlowVar}{}\)
given in section \ref{ssec::MeanOp} we obtain
that for each \(\FlowVar>0\)
\begin{equation*}
	\DeltaR{\VFieldE}{\FlowVar}{\funcA}=
	\lim_{\pFVar\to0}\MeanOp{\VFieldE}{\pFVar}{\DeltaR{\VFieldE}{\FlowVar}{\funcA}}=
	\lim_{\pFVar\to0}\MeanOp{\VFieldE}{\FlowVar}{\DeltaR{\VFieldE}{\pFVar}{\funcA}}=
	\MeanOp{\VFieldE}{\FlowVar}{\funcB},
\end{equation*}
that is \(\funcB=L_\VFieldE\funcA\).

%\Case{\stmGFL}\(\Longrightarrow\)\Case{\stmGFZ}:
%the equality \(\funcB=L_\VFieldE\funcA\) may be rephrased saying that
%when \(\FlowVar\neq0\) then
%\begin{equation*}
%	\DeltaR{\VField}{\FlowVar}{\funcA}=\MeanOp{\VField}{\FlowVar}{\funcB}.
%\end{equation*}
%Then Proposition \ref{stm::MeanOpMain} implies that
%\begin{equation*}
%	\lim_{\FlowVar\to0}\DeltaR{\VField}{\FlowVar}{\funcA}=
%	\lim_{\FlowVar\to0}\MeanOp{\VField}{\FlowVar}{\funcB}=
%	\funcB,
%\end{equation*}
%that is \(\funcA\in\Dom(\GFZ)\) and \(\GFZ\funcA=\funcB\).

\Case{\stmGFL}\(\Longrightarrow\)\Case{\stmGFX}:
if
\(\funcB=L_\VFieldE\funcA\)
then Proposition \ref{stm::AlmostLieToLie} implies that
\(\DeltaR{\VFieldE}{\FlowVar}{\funcA}=\MeanOp{\VFieldE}{\FlowVar}{\funcB}\)
for each \(\FlowVar>0\).

Let \(\funcTest\in\proofTestSpace\) be given.
Propositions \ref{stm::DQLim} and \ref{stm::MeanOpMain} imply respectively
\begin{equation*}
	\lim_{\FlowVar\to0}\int_\OpenE\bigl(\DeltaR\VField\FlowVar\funcA\bigr)\funcTest=
	-\int_\OpenE\funcA\VFieldE\funcTest
	-\int_\OpenE\funcA\funcTest\Divergence\VFieldE
\end{equation*}
and
\begin{equation*}
	\lim_{\FlowVar\to0}\int_\OpenE\bigl(\MeanOp{\VFieldE}{\FlowVar}\funcB\bigr)\funcTest
	=\int_\OpenE\funcB\funcTest,
\end{equation*}
and hence
\begin{equation*}
	\int_\OpenE\funcB\funcTest=
	-\int_\OpenE\funcA\VFieldE\funcTest
	-\int_\OpenE\funcA\funcTest\Divergence\VFieldE.
\end{equation*}
Since \(\funcTest\in\proofTestSpace\) is arbitrary
it follows that
\(\funcB=\VField\funcA\).

}\end{proof}

}%end unit
%}{}

\subsection{Existence of the derivative}
% !TeX encoding = UTF-8
% !TeX spellcheck = en_GB
{%begin unit
\def\FVarUpper{{\FlowVar_0}}
Now we show that some compactness properties of the different quotients
imply the existence of the derivatives.
%Let \(\VField\)  be a complete Lipshitzian vector field
%on an open set \(\OpenE\subset\RR^\Dim\).
%Let \(\GPFlow{\VFieldE}{\FlowVar}\) be the induced
%\(c_0\) one parameter group on the Banach space \(L^1(\OpenE)\),
%\begin{equation*}
%	L^1(\OpenE)\ni\funcA\mapsto\GPFlow\VField\FlowVar\funcA=\funcA\circ\VFlowW{\VField}{\FlowVar}\in L^1(\OpenE)
%\end{equation*}
%and let \(\GFZ:\Dom(\GFZ)\to L^1(\OpenE)\) be its infinitesimal
%generator.

\begin{theorem}\label{stm::DerExistDQComp}
A function \(\funcA\in L^1(\OpenE)\) admits the Lie/distributional
derivative in \(L^1(\OpenE)\) if, and only if, for some \(\FVarUpper>0\)
the family of the different quotients
\begin{equation*}
	\DeltaR{\VFieldE}{\FlowVar}{\funcA}
	=\dfrac{\funcA\bigl(\VFlow{\VField}{\FlowVar}{\cdot}\bigr)-\funcA(\cdot)}{\FlowVar}
	,\ 0<\FlowVar<\FVarUpper
\end{equation*}
is weakly relatively compact in \(L^1(\OpenE)\).
\end{theorem}

\begin{proof}{
\def\aLeft{a}
\def\bRight{b}
\def\Idx{k}
Assume that the function \(\funcA\) admits the Lie derivative and
\(g=L_\VFieldE\funcA\). Proposition \ref{stm::AlmostLieToLie} implies that for each \(\FlowVar\in\RR\) we have
\begin{equation*}
	\DeltaR{\VField}{\FlowVar}{\funcA}=\MeanOp{\VField}{\FlowVar}{\funcB}
\end{equation*}
and by Proposition \ref{stm::MeanOpMain} the map
\begin{equation*}
	\RR\ni\FlowVar\mapsto\MeanOp{\VFieldE}{\FlowVar}\funcA\in L^1(\OpenE)
\end{equation*}
is continuous.
Since the image of any compact set under a continuous map is compact,
it follows that for each \(-\infty<\aLeft<\bRight<+\infty\) the family
\begin{equation*}
	\DeltaR{\VFieldE}{\FlowVar}{\funcA}
	,\ \aLeft<\FlowVar<\bRight
\end{equation*}
is relatively compact in \(L^1(\OpenE)\)
and hence
is weakly relatively compact in \(L^1(\OpenE)\).

Conversely, assume that the family
\begin{equation*}
	\DeltaR{\VFieldE}{\FlowVar}{\funcA}
	,\ 0<\FlowVar<\FVarUpper
\end{equation*}
is weakly relatively compact in \(L^1(\OpenE)\) for some \(\FVarUpper>0\).

The classical Eberlein-Smulian theorem implies that such a family also is
weakly sequentially compact
and hence it is possible to find 
a sequence \(\FlowVar_\Idx>0\) converging to \(0\)
such that
%the sequence of
the difference quotients
\(\DeltaR{\VFieldE}{\FlowVar_\Idx}{\funcA}\)
converges weakly in \(L^1(\OpenE)\) when \(\Idx\to+\infty\)
to a function \(\funcB\in L^1(\OpenE)\).
%there exist a sequence \(\FlowVar_\Idx>0\) which converges to \(0\)
%and
%\begin{equation*}
%	\funcB=*-\lim_{\Idx\to\infty}\DeltaR{\VFieldE}{\FlowVar_\Idx}{\funcA}
%\end{equation*}
%for some \(\funcB\in L^1(\OpenE)\).
Then Proposition \ref{stm::GFWStar} immediately yields that
\(\funcA\in\Dom(\GFZ)\) and \(\funcB=\GFZ\funcA\).
Theorem \ref{stm::LOne::ZEqLEqX} concludes the proof.

}\end{proof}
}%end unit
\subsection{Upper gradients}
% !TeX encoding = UTF-8
% !TeX spellcheck = en_GB
{%begin unit
\def\FFamA{\mathcal{F}}
\def\FFamB{\mathcal{G}}
%Inspired by \cite{mem:HajlaszKoskela:SobolevMetPoincare:MR1683160}
%we say that a non negative function \(\UGradA\in\LOneLoc{\OpenE}\)
%is an \emph{upper gradient} for the function \(\funcA\in\LOneLoc{\OpenE}\)
%along the vector field \(\VFieldE\)
%if for each \(\FlowVar>0\)
%\begin{equation*}
%	\abs{\funcA\bigl(\VFlow{\VField}{\FlowVar}{\pEuc}\bigr)-\funcA(\pEuc)}
%	\leq
%	\int_0^\FlowVar\UGradA\bigl(\VFlow{\VFieldE}{\FlowVarB}{\pEuc}\bigr)\,d\FlowVarB
%\end{equation*}
Now we prove the relationships between derivatives
and upper gradients along the vector field \(\VField\).

Observe that Proposition \ref{stm::AlmostUGToUG} implies that
\(\UGradA\in\LOneLoc{\OpenE}\)
is an upper gradient for the function \(\funcA\in\LOneLoc{\OpenE}\)
along the vector field \(\VFieldE\)
if, and only if,
for each \(\FlowVar>0\)
%for almost all \(\pEuc\in\OpenE\), i. e.
\begin{equation*}
	\abs{\DeltaR{\VFieldE}{\FlowVar}{\funcA}}
	\leq
	\MeanOp{\VFieldE}{\FlowVar}{\UGradA}
\end{equation*}
almost everywhere in \(\OpenE\).

\begin{lemma}{
\def\fa{u}
\def\fb{v}
Let \(\FFamA,\FFamB\subset L^1(\OpenE)\) be two families of functions.
If the family \(\FFamB\) is weakly compact in \(L^1(\OpenE)\)
and
for each \(\fa\in\FFamA\) there exists \(\fb\in\FFamB\) such that
\begin{equation*}\
	\abs{\fa}\leq\abs{\fb}
\end{equation*}
almost everywhere on \(\OpenE\)
then the family \(\FFamA\) is also weakly compact in \(L^1(\OpenE)\).
}\end{lemma}

\begin{proof}
It is immediate to show that the family \(\FFamB\) satisfies
the conditions of the Dunford-Pettis theorem
(Theorem \ref{stm::DunfordPettis}).

\end{proof}

\begin{theorem}\label{stm::LOne::UGradToX}
If the function \(\funcA\in L^1(\OpenE)\)
admits an upper gradient \(\UGradA\in L^1(\OpenE)\)
along the vector field \(\VFieldE\)
then it also admits the distributional derivative
\(\funcB=\VFieldE\funcA\in L^1(\OpenE)\).
%with respect to \(\VFieldE\).
\end{theorem}

\begin{proof}{
\def\tz{\FlowVar_0}
Let us fix \(\tz>0\).
We have already observed in the proof of Theorem \ref{stm::DerExistDQComp}
that Proposition \ref{stm::MeanOpMain}
implies the weakly compactness in \(L^1(\OpenE)\) of the family
\begin{equation*}
	\MeanOp{\VFieldE}{\FlowVar}{\UGradA},\ 0<\FlowVar<\tz.
\end{equation*}
By hypotheses, we have
\begin{equation*}
	\abs{\DeltaR{\VFieldE}{\FlowVar}{\funcA}}
	%\abs{\funcA\bigl(\VFlow{\VField}{\FlowVar}{\pEuc}\bigr)-\funcA(\pEuc)}
	\leq
	\MeanOp{\VFieldE}{\FlowVar}{\UGradA},\ 0<\FlowVar<\tz
\end{equation*}
almost everywhere on \(\OpenE\) and hence,
by the previous lemma the family
\begin{equation*}
	\DeltaR{\VFieldE}{\FlowVar}{\funcA}
	%\funcA\bigl(\VFlow{\VField}{\FlowVar}{\pEuc}\bigr)-\funcA(\pEuc)
	,\ 0<\FlowVar<\tz
\end{equation*}
is weakly compact in \(L^1(\OpenE)\).

Theorem \ref{stm::DerExistDQComp}
completes the proof.

}\end{proof}

Conversely we have:

\begin{theorem}\label{stm::LOne::XToUGrad}
Assume that the function \(\funcA\in L^1(\OpenE)\)
admits the distributional derivative
\(\funcB=\VFieldE\funcA\in L^1(\OpenE)\)
with respect to \(\VFieldE\).
Then \(\abs{\funcB}\)
is an upper gradient of \(\funcA\)
along the vector field \(\VFieldE\)
and every upper gradient \(\UGradA\in L^1(\OpenE)\)
of \(\funcA\)
along the vector field \(\VFieldE\)
satisfies
\begin{equation*}
	\abs{\funcB}\leq\UGradA
\end{equation*}
almost everywhere on \(\OpenE\).
\end{theorem}

\begin{proof}
Since \(\funcB=\VFieldE\funcA\)
then, by Theorem \ref{stm::LOne::ZEqLEqX},
we also have \(\funcB=L_\VFieldE\funcA\),
that is
\begin{equation*}
	\DeltaR{\VField}{\FlowVar}{\funcA}(\pEuc)=\MeanOp{\VField}{\FlowVar}{\funcB}(\pEuc)
\end{equation*}
%when \(\FlowVar>0\)
and hence
\begin{equation*}
	\abs{\DeltaR{\VField}{\FlowVar}{\funcA}(\pEuc)}
	=
	\abs{\MeanOp{\VField}{\FlowVar}{\funcB}(\pEuc)}
	\leq
	\MeanOp{\VField}{\FlowVar}{\abs{\funcB}(\pEuc)}	
\end{equation*}
for almost all \((\FlowVar,\pEuc)\in\DFlowX\VField\),
that is \(\abs{\funcB}\)
is an upper gradient of \(\funcA\).

Suppose now that \(\UGradA\in L^1(\OpenE)\)
is an arbitrary upper gradient of \(\funcA\)
along the vector field \(\VFieldE\).

Then when \(\FlowVar>0\) we have by definition
\begin{equation*}
	\abs{\DeltaR{\VFieldE}{\FlowVar}{\funcA}}
	\leq
	\MeanOp{\VFieldE}{\FlowVar}{\UGradA}
\end{equation*}
almost everywhere on \(\OpenE\).
Taking the limit as \(\FlowVar\to0^+\),
recalling that
\begin{equation*}
	\funcB=L_\VFieldE\funcA
	=
	\lim_{\FlowVar\to0^+}\DeltaR{\VField}{\FlowVar}{\funcA}
\end{equation*}
%with respect to the topology in \(L^1(\OpenE)\) induced by its standard norm,
we obtain
\begin{equation*}
	\abs{\funcB}
	=
	\lim_{\FlowVar\to0^+}\abs{\DeltaR{\VField}{\FlowVar}{\funcA}}
	\leq
	\lim_{\FlowVar\to0^+}\MeanOp{\VField}{\FlowVar}{\UGradA}
	=\UGradA
\end{equation*}
almost everywhere on \(\OpenE\),
the last equality being the consequence of Proposition \ref{stm::MeanOpMain}.

\end{proof}

}%end unit
%\ifdef{\cbFullBase}{%}{}
%\subsection{Meyer-Serrin theorems}
%\input{src/VFLieDerMS}
%}{}

\section{\label{section:MainProof}Proofs of main theorems}
\subsection{Localization lemmas}
% !TeX encoding = UTF-8
% !TeX spellcheck = en_GB

{%begin unit
\def\funcAA{{\funcA_1}}

Let \(\VField\) be a locally Lipshitzian vector field on the open set \(\OpenE\subset\RR^\Dim\)
with associated local flow %\(\VFlow{\VFieldE}{\FlowVar}{\pEuc}\)
\begin{equation*}
	\RR\times\OpenE\supset\DFlowX\VFieldE\ni(\FlowVar,\pEuc)\mapsto\VFlow{\VFieldE}{\FlowVar}{\pEuc}\in\OpenE
\end{equation*}
and let \(\funcA\in\LOneLoc\OpenE\) be given.
We need the following lemmas.

\begin{lemma}
Let \(\CutOffA\in\LipC{\OpenE}\)
be a nonnegative function
and suppose that
\(L_\VFieldE\funcA\in\LOneLoc{\OpenE}\).
%\(\VFieldB=\CutOffA\VFieldE\)
Then
\begin{equation*}
	L_{\CutOffA\VFieldA}\funcA=\CutOffA L_\VFieldE\funcA\in\LOneLoc{\OpenE}.
\end{equation*}
\end{lemma}

\begin{proof}
{

The vector field \(\VFieldB=\CutOffA\VFieldE\)
has compact support in \(\OpenE\)
and hence is complete.

%(see the next session)
It is not difficult to prove
that the non negativity and the compactness of the support of \(\CutOffA\)
imply that there exists
a unique function
\begin{equation*}
	\ChParam:\RR\times\OpenE\to\RR
\end{equation*}
which verifies
\begin{equation*}
	\begin{cases}
		\dfrac{\partial\ChParam(\xVar,\pEuc)}{\partial\xVar}
		=\Cutoff\bigl(\VFlow{\VFieldE}{\ChParam(\xVar,\pEuc)}{\pEuc}\bigr)
		\\
		\ChParam(0,\pEuc)=0
	\end{cases}
	.
\end{equation*}
%For convenience the proof of this fact will be given in the next section.

The flow \(\VFlow{\VFieldB}{\FlowVar}{\pEuc}\)
associated to \(\VFieldB\)
is then given by
\begin{equation*}
	\VFlow{\VFieldB}{\FlowVar}{\pEuc}
	=
	\VFlow{\VFieldE}{\ChParam(\FlowVar,\pEuc)}{\pEuc}.
\end{equation*}
Since
\begin{equation*}
	\Cutoff(\pEuc)>0\ \Longrightarrow \dfrac{\partial\ChParam(\xVar,\pEuc)}{\partial\xVar}>0,\ \xVar\in\RR
\end{equation*}
and
\begin{equation*}
	\Cutoff(\pEuc)=0\ \Longrightarrow \ChParam(\xVar,\pEuc)=0,\ \xVar\in\RR,
\end{equation*}
it follows that
for almost all \((\FlowVar,\pEuc)\in]0,+\infty[\times\OpenE\)
we have
\begin{equation*}
	\begin{split}
		\funcA\bigl(\VFlow{\VFieldB}{\FlowVar}{\pEuc}\bigr)
		-\funcA(\pEuc)
		&=
		\funcA\bigl(\VFlow{\VField}{\ChParam(\FlowVar,\pEuc)}{\pEuc}\bigr)
		-\funcA(\pEuc)
		%\\
		=\int_{0}^{\ChParam(\FlowVar,\pEuc)}
		L_\VFieldE\funcA\bigl(\VFlow{\VFieldE}{\FlowVarB}{\pEuc}\bigr)\,d\FlowVarB
		\\
		&
		=\int_{0}^{\FlowVar}
		L_\VFieldE\funcA\bigl(\VFlow{\VFieldE}{\ChParam(\FlowVarB,\pEuc)}{\pEuc}\bigr)
		\dfrac{\partial\ChParam(\FlowVarB,\pEuc)}{\partial\FlowVarB}
		\,d\FlowVarB
		\\
		&
		=\int_{0}^{\FlowVar}
		\Cutoff\bigl(\VFlow{\VFieldE}{\ChParam(\FlowVarB,\pEuc)}{\pEuc}\bigr)
		L_\VFieldE\funcA\bigl(\VFlow{\VFieldE}{\ChParam(\FlowVarB,\pEuc)}{\pEuc}\bigr)
		\,d\FlowVarB
		\\
		&
		=\int_{0}^{\FlowVar}
		\Cutoff\bigl(\VFlow{\VFieldB}{\FlowVarB}{\pEuc}\bigr)
		L_\VFieldE\funcA\bigl(\VFlow{\VFieldB}{\FlowVarB}{\pEuc}\bigr)
		\,d\FlowVarB.
	\end{split}
\end{equation*}
and hence \(L_\VFieldB\funcA=\CutOffA L_\VFieldE\funcA\),

}

\end{proof}

\begin{lemma}
Let \(\CutOffA\in C^\infty(\OpenE)\)
be given
and suppose that
\(\VFieldE\funcA\in\LOneLoc{\OpenE}\).
%\(\VFieldB=\CutOffA\VFieldE\)
Then
\begin{equation*}
	(\CutOffA\VFieldE)\funcA=\CutOffA\VFieldE\funcA\in\LOneLoc{\OpenE}.
\end{equation*}
\end{lemma}

\begin{proof}		
Let \(\TestFEuc\in C_0^\infty(\OpenE)\) be given. We have
\begin{equation*}
	\begin{split}
		\Int(\CutOffA\VFieldE\funcA)\TestFEuc
		&=
		\Int\VFieldE\funcA(\CutOffA\TestFEuc)
		=
		-\Int\funcA\VField(\CutOffA\TestFEuc)
		-\Int\funcA\Cutoff\TestFEuc\Divergence\VFieldE
		\\
		&=
		-\Int\funcA(\VField\CutOffA)\TestFEuc
		-\Int\CutOffA\funcA(\VField\TestFEuc)
		-\Int\funcA\Cutoff\TestFEuc\Divergence\VFieldE
		\\
		&=
		-\Int\CutOffA\funcA(\VField\TestFEuc)
		-\Int\funcA(\VField\CutOffA+\Cutoff\Divergence\VFieldE)\TestFEuc
		\\
		&=
		-\Int\funcA(\CutOffA\VField\TestFEuc)
		-\Int\funcA\Divergence(\Cutoff\VFieldE)\TestFEuc.
	\end{split}
\end{equation*}
Since \(\TestFEuc\in C_0^\infty(\OpenE)\) is arbitrary,
the assertion follows.
		
\end{proof}

%The proof of the following lemma is straightforward and we omit it.

\begin{lemma}
If \(\funcAA\in\LOneLoc{\OpenE}\) satisfies
\(\funcAA(\pEuc)=\funcA(\pEuc)\) whenever \(\VFieldE(\pEuc)\neq0\)
then:

%\begin{tabular}{rl}
%	\Case{1}&if \(L_\VField\funcA\in\LOneLoc{\OpenE}\) then \(L_\VFieldE\funcA_1=L_\VFieldE\funcA\);\\ 
%	\Case{2}&if \(\VField\funcA\in\LOneLoc{\OpenE}\) then \(\VFieldE\funcA_1=\VFieldE\funcA\).
%\end{tabular}
%
\begin{enumerate}[label=\cbRomanLabel]
	\item\label{stm::LocalA}
	if \(L_\VField\funcA\in\LOneLoc{\OpenE}\) then \(L_\VFieldE\funcAA=L_\VFieldE\funcA\);
	\item\label{stm::LocalB}
	if \(\VField\funcA\in\LOneLoc{\OpenE}\) then \(\VFieldE\funcAA=\VFieldE\funcA\).
\end{enumerate}
\end{lemma}

\begin{proof}{
\def\NZetaRegion{U}
The set
\begin{equation*}
	\NZetaRegion=\bigl\{\pEuc\in\OpenE\mid\VFieldE(\pEuc)\neq0\bigr\}
\end{equation*}
is an open subset of \(\OpenE\) which is invariant under the
flow \(\VFlow{\VFieldE}{\FlowVar}{\pEuc}\)
and hence \(L_\VFieldE\funcA_1(\pEuc)=L_\VFieldE\funcA(\pEuc)\)
if \(\pEuc\in\NZetaRegion\).

Observe that \(\VFlow{\VFieldE}{\FlowVar}{\pEuc}=\pEuc\)
whenever \(\pEuc\in\OpenE\setminus\NZetaRegion\) and hence
both \(L_\VFieldE\funcA\) and \(L_\VFieldE\funcAA\) vanish on \(\pEuc\in\OpenE\setminus\NZetaRegion\).

This conclude the proof of \ref{stm::LocalA}.

Using the coarea formula for Lipshitzian maps
(see, e.g., \cite[Theorem 3.2.11]{book:Federer})
it is not difficult to prove that all the (partial) derivatives
of the components of the vector field \(\VFieldE\)
vanish almost everywhere on \(\OpenE\setminus\NZetaRegion\) and hence
\begin{equation*}
	\Divergence\VFieldE(\pEuc)=0
\end{equation*}
for almost \(\pEuc\in\OpenE\setminus\NZetaRegion\).
It follows that for each test function
\(\TestFEuc\in C_0^\infty(\OpenE)\)
\begin{equation*}
	\begin{split}
		\Int(\VFieldE\funcA)\TestFEuc
		&=
		-\Int\funcA\VFieldE\TestFEuc-\Int\funcA\TestFEuc\Divergence\VFieldE
		=
		-\int_\NZetaRegion\funcA\VFieldE\TestFEuc-\int_\NZetaRegion\funcA\TestFEuc\Divergence\VFieldE
		\\
		&=
		-\int_\NZetaRegion\funcAA\VFieldE\TestFEuc-\int_\NZetaRegion\funcAA\TestFEuc\Divergence\VFieldE
		\\
		&=
		-\Int\funcAA\VFieldE\TestFEuc-\Int\funcAA\TestFEuc\Divergence\VFieldE	
	\end{split}
\end{equation*}
and hence \ref{stm::LocalB} follows.
}

\end{proof}

}%end unit
\subsection{The proofs}
%\ifdef{\cbVOne}{
%	\input{src/VFMainProof}
%}{}
%\ifdef{\cbVTwo}{
%	\input{src/VFMainProofVTwo}
%}{}
%\ifdef{\cbVThree}{
% !TeX encoding = UTF-8
% !TeX spellcheck = en_GB

{%begin unit
\def\funcX{\tilde{\funcA}}
\def\CoeffBall{S^{\VFIMax}}
\def\SetEstimateOk{E}
\def\OpenEOk{\tilde{\OpenE}}
\def\vEucA{v}
\def\vEucB{w}
\def\kIdx{k}

%Be begin now proving Theorem \ref{stm::MainThm}.

We are ready to prove Theorem \ref{stm::MainThm}.

We will prove that \Case{\stmGFXLoc}\(\Longrightarrow\)\Case{\stmGFLLoc}\(\Longrightarrow\)\Case{\stmGFZLoc}\(\Longrightarrow\)\Case{\stmGFXLoc}.

\Case{\stmGFXLoc}\(\Longrightarrow\)\Case{\stmGFLLoc}:
assume that \(\funcB=\VFieldE\funcA\in\LOneLoc{\OpenE}\).
We then prove that \(\funcB\) coincides with
the Lie derivative of \(\funcA\) with respect to \(\VFieldE\).

It is possible to find a sequence \(\FlowVar_\kIdx\) of positive numbers and
a sequence \(\SubOpenE_\kIdx\) of open relatively compact subset of \(\OpenE\)
such that
\begin{equation*}
	\DFlowXPlus\VFieldE=\bigcup_{\kIdx=1}^\infty]0,\FlowVar_\kIdx[\times\SubOpenE_\kIdx
\end{equation*}
and \([0,\FlowVar_0]\times\overline\SubOpenE\subset\DFlowX\VFieldE\) for each \(\kIdx\).

It suffices then to prove that
\begin{equation*}
	\funcA\bigl(\VFlow{\VFieldE}{\FlowVar}{\pEuc}\bigr)-\funcA(\pEuc)
	=\int_0^\FlowVar
	\funcB\bigl(\VFlow{\VFieldE}{\FlowVarB}{\pEuc}\bigr)\,d\FlowVarB
\end{equation*}
for almost all \((\FlowVar,\pEuc)\in]0,\FlowVar_\kIdx[\times\SubOpenE_\kIdx\)
for each fixed \(\kIdx\).

%Since the open set \(\SubOpenE\) can be
%covered by a countable family of relatively compact open subset of \(\OpenE\)
%we may assume that \(\SubOpenE\) itself is relatively compact in \(\OpenE\)
%and \([0,\FlowVar_0]\times\overline\SubOpenE\subset\DFlowX\VFieldE\).
%
So let \(\kIdx\) be fixed and
let \(\CutOffA\in C_0^\infty(\OpenE)\) be a function
which satisfies \(\Cutoff\equiv1\) in a neighbourhood of 
the relatively compact set
\begin{equation*}
	\bigcup_{\FlowVar\in[0,\FlowVar_0]}\VFlow{\VFieldE}{\FlowVar}{\SubOpenE_\kIdx}
\end{equation*}
and
let \(\funcX\) the function defined by
\begin{equation*}
	\funcX(\pEuc)=
	\begin{cases}
		\funcA(\pEuc)\ &\text{if}\ \CutOffA(x)>0;\\
		0\ &\text{if}\ \CutOffA(x)=0.\\
	\end{cases}
\end{equation*}

By the Lemmas of the previous subsection we obtain
\begin{equation*}
	\CutOffA\funcB
	=\CutOffA(\VFieldE\funcA)
	=(\CutOffA\VFieldE)\funcA
	=(\CutOffA\VFieldE)\funcX.
\end{equation*}
Since \(\CutOffA\VFieldE\) is complete and \(\funcX, \CutOffA\funcB\in L^1(\OpenE)\),
then Theorem \ref{stm::LOne::ZEqLEqX} implies that
\begin{equation*}
	\CutOffA\funcB
	=(\CutOffA\VFieldE)\funcX
	=L_{\CutOffA\VFieldE}\funcX.
\end{equation*}
Since \(\CutOffA\bigl(\VFlow{\VFieldE}{\VFieldE}{\pEuc}\bigr)=1\)
for each \(\FlowVar\in[0,\FlowVar_\kIdx]\) and each \(\pEuc\in\SubOpenE\)
it follows that
\begin{equation*}
	[0,\FlowVar_\kIdx]\ni\FlowVar\mapsto\VFlow{\VFieldE}{\FlowVar}{\pEuc}
\end{equation*}
is an integral curve of the vector field \(\CutOffA\VFieldE\)
and hence
\begin{equation*}
	\begin{split}
		\funcA\bigl(\VFlow{\VFieldE}{\FlowVar}{\pEuc}\bigr)
		-\funcA(\pEuc)
		&=\funcX\bigl(\VFlow{\CutOffA\VFieldE}{\FlowVar}{\pEuc}\bigr)
		-\funcX(\pEuc)
		=\int_0^\FlowVar
		\CutOffA\bigl(\VFlow{\CutOffA\VFieldE}{\FlowVarB}{\pEuc}\bigr)
		\funcB\bigl(\VFlow{\CutOffA\VFieldE}{\FlowVarB}{\pEuc}\bigr)\,d\FlowVarB
		\\
		&
		=\int_0^\FlowVar
		\funcB\bigl(\VFlow{\VFieldE}{\FlowVarB}{\pEuc}\bigr)\,d\FlowVarB
	\end{split}
\end{equation*}
for almost all \((\FlowVar,\pEuc)\in]0,\FlowVar_\kIdx[\times\SubOpenE_\kIdx\),
as required.

\Case{\stmGFLLoc}\(\Longrightarrow\)\Case{\stmGFZLoc}:
assume that \(\funcB=L_\VFieldE\funcA\)
and let \(\SubOpenE\) be a relatively compact subset of \(\OpenE\).
Then we have \([0,\FlowVar_0]\times\overline\SubOpenE\subset\DFlowX\VFieldE\)
for some \(\FlowVar_0>0\).
We then obtain
\begin{equation*}
	\CutOffA\funcB
	=\CutOffA(L_\VFieldE\funcA)
	=L_{\CutOffA\VFieldE}\funcA
	=L_{\CutOffA\VFieldE}\funcX,
\end{equation*}
where \(\Cutoff\) and \(\funcX\) are defined as before.

Theorem \ref{stm::LOne::ZEqLEqX} again implies that
\(%\begin{equation*}%\label{eq::LieAlmostEverywhere}
\lim_{\FlowVar\to0^+}
\DeltaR{\CutOffA\VField}\FlowVar\funcX=\CutOffA\funcB
%\dfrac{\funcA\circ\VFlowMap{\VField}{\FlowVar}-\funcA}{\FlowVar}=\funcB
%\FlowVar^{-1}\bigl(\funcA\circ\VFlowMap{\VField}{\FlowVar}-\funcA\bigr)=\funcB
\) %\end{equation*}
with respect to the \(L^1(\OpenE)\) topology.
Observing that
\(\CutOffA\bigl(\VFlow{\VFieldE}{\FlowVar}{\pEuc}\bigr)=1\)
for each \(\FlowVar\in[0,\FlowVar_0]\) and each \(\pEuc\in\SubOpenE\)
it follows that
\begin{equation*}%\label{eq::LieAlmostEverywhere}
	\lim_{\FlowVar\to0^+}\int_{\SubOpenE}
	\abs{\dfrac{\funcA\bigl(\VFlow{\VField}{\FlowVar}{\pEuc}\bigr)-\funcA(\pEuc)}{\FlowVar}-\funcB(\pEuc)}\,d\pEuc=0,
\end{equation*}
as desired.

\Case{\stmGFZLoc}\(\Longrightarrow\)\Case{\stmGFXLoc}:
follows immediately from Proposition \ref{stm::DQLim}.

%The validity of the formula \eqref{eq::LieAlmostEverywhere}
%follows easily from Theorem \ref{stm::LOne::ZEqLEqX}.
%see also Proposition \ref{stm::LIpVF::PointwiseLieDer}.

The proof of
Theorems \ref{stm::Main::UGradToX} and \ref{stm::Main::XToUGrad}
can be obtained by similar arguments
using respectively
Theorem \ref{stm::LOne::UGradToX}
and
Theorem \ref{stm::LOne::XToUGrad}.

Now we prove Theorem \ref{stm::Main::UGradToXSystem}.
Let
\begin{equation*}
	\VFieldE_\VFIMin,\ldots,\VFieldE_\VFIMax,
\end{equation*}
and \(\UGradA,\funcA\in\LOneLoc{\OpenE}\)
be as in the theorem. % \ref{stm::Main::UGradToXSystem}.

Let \(\VConstX_\VFIMin,\ldots,\VConstX_\VFIMax\) be given satisfying
\begin{equation*}
	\sum_{\VFIndex=\VFIMin}^{\VFIMax}\VConstX_\VFIndex^2\leq1.
\end{equation*}

Theorem \ref{stm::Main::UGradToX} implies that
the distributional derivatives
\(\VFieldE_\VFIMin\funcA,\ldots,\VFieldE_\VFIMax\funcA\)
and
\(\Bigl(\sum_{\VFIndex=\VFIMin}^{\VFIMax}\VConstX_\VFIndex\VFieldE_\VFIndex\Bigr)\funcA\)
are locally integrable and
\begin{equation*}
	\abs{\Bigl(\sum_{\VFIndex=\VFIMin}^{\VFIMax}\VConstX_\VFIndex\VFieldE_\VFIndex\Bigr)\funcA}\leq\UGradA
\end{equation*}
almost everywhere on \(\OpenE\).

Of course we have
\begin{equation*}
	\Bigl(\sum_{\VFIndex=\VFIMin}^{\VFIMax}\VConstX_\VFIndex\VFieldE_\VFIndex\Bigr)\funcA
	=
	\sum_{\VFIndex=\VFIMin}^{\VFIMax}\VConstX_\VFIndex\VFieldE_\VFIndex\funcA
\end{equation*}
as distributional derivatives and hence as \(\LOneLoc{\OpenE}\) functions
and therefore we obtain
\begin{equation*}
	\abs{\sum_{\VFIndex=\VFIMin}^{\VFIMax}\VConstX_\VFIndex\VFieldE_\VFIndex\funcA}
	\leq\UGradA
\end{equation*}
almost every where on \(\OpenE\).

Let now
\begin{equation*}
	\CoeffBall=\bigl\{\VConstX\in\RR^\VFIMax\mid
	\sum_{\VFIndex=\VFIMin}^{\VFIMax}\VConstX_\VFIndex^2\leq1\bigr\}
\end{equation*}
and let \(\SetEstimateOk\subset\CoeffBall\times\OpenE\)
be the set of the pairs \((\VConstX,\pEuc)\) where the inequality
\begin{equation}\label{eq::ineqA}
	\abs{\sum_{\VFIndex=\VFIMin}^{\VFIMax}\VConstX_\VFIndex\VFieldE_\VFIndex\funcA(\pEuc)}
	\leq\UGradA(\pEuc)
\end{equation}
is satisfied.
Clearly \(\SetEstimateOk\) is a measurable subset of \(\CoeffBall\times\OpenE\)
and we have shown that for each \(\VConstX\in\CoeffBall\)
the set of \(\pEuc\in\OpenE\) such that \((\VConstX,\pEuc)\in\SetEstimateOk\)
has full measure in \(\OpenE\).

By the Fubini-Tonelli theorems for almost all \(\pEuc\in\OpenE\)
the inequality \eqref{eq::ineqA} is satisfied for almost all \(\VConstX\in\CoeffBall\).

Observing that for each \(\pEuc\in\OpenE\) fixed the function
\begin{equation*}
	\CoeffBall\ni\VConstX\mapsto
	\abs{\sum_{\VFIndex=\VFIMin}^{\VFIMax}\VConstX_\VFIndex\VFieldE_\VFIndex\funcA(\pEuc)}
	\in\RR
\end{equation*}
is continuous we can conclude that there is
a subset \(\OpenEOk\subset\OpenE\) of full measure in \(\OpenE\)
such that for all \(\pEuc\in\OpenEOk\) the inequality \eqref{eq::ineqA}
is satisfied for all \(\VConstX\in\CoeffBall\).

Using the standard equality
\begin{equation*}
	\Bigl(\sum_{\VFIndex=\VFIMin}^{\VFIMax}\vEucA_\VFIndex^2\Bigr)^{1/2}=
	\sup_{\vEucB\in\CoeffBall}\sum_{\VFIndex=\VFIMin}^{\VFIMax}\vEucA_\VFIndex\vEucB_\VFIndex
\end{equation*}
we obtain that when \(\pEuc\in\OpenEOk\)
\begin{equation*}
	\Bigl(\sum_{\VFIndex=\VFIMin}^{\VFIMax}\VFieldE_\VFIndex\funcA(\pEuc)^2\Bigr)^{1/2}=
	\sup_{\VConstX\in\CoeffBall}\sum_{\VFIndex=\VFIMin}^{\VFIMax}\VConstX_\VFIndex\VFieldE_\VFIndex\funcA(\pEuc)
	\leq\UGradA(\pEuc),
\end{equation*}
as desired.

}%end unit
%}{}
%\ifdef{\cbAppendixA}{
%	\subsection{A Cauchy problem}
%	\input{src/AppendixCompleteX}
%}{}

\bibliographystyle{amsalpha}

\providecommand{\bysame}{\leavevmode\hbox to3em{\hrulefill}\thinspace}
\providecommand{\MR}{\relax\ifhmode\unskip\space\fi MR }
\providecommand{\MRhref}[2]{%
  \href{http://www.ams.org/mathscinet-getitem?mr=#1}{#2}
}
\providecommand{\href}[2]{#2}

\end{document}